\documentclass[10pt,letterpaper]{article}
\usepackage[top=0.85in,left=2.75in,footskip=0.75in]{geometry}

\usepackage{multicol}

\usepackage{booktabs}
\usepackage{longtable}

\usepackage{tikz}
\usetikzlibrary{arrows,shapes,quotes}
\usetikzlibrary{patterns,decorations.pathreplacing}

\usepackage{algorithm}
\usepackage{algpseudocode}

\usepackage{placeins}
\usepackage{float}

\usepackage{amsmath,amssymb}
\allowdisplaybreaks

\usepackage{changepage}

\usepackage{textcomp,marvosym}

\usepackage{cite}

\usepackage{nameref,hyperref}

\usepackage{mathtools}
\usepackage{blkarray, bigstrut}
\usepackage{multirow}
\usepackage{array}

\usepackage[right]{lineno}

\usepackage[nopatch=eqnum]{microtype}
\DisableLigatures[f]{encoding = *, family = * }

\usepackage[table]{xcolor}

\usepackage{array}

\newcolumntype{+}{!{\vrule width 2pt}}

\newlength\savedwidth

\raggedright
\usepackage[aboveskip=1pt,labelfont=bf,labelsep=period,justification=raggedright,singlelinecheck=off]{caption}

\makeatletter
\renewcommand{\@biblabel}[1]{\quad#1.}
\makeatother

\usepackage{lastpage,fancyhdr,graphicx}
\usepackage{epstopdf}
\fancyheadoffset[L]{2.25in}
\fancyfootoffset[L]{2.25in}
\newtheorem{theorem}{Theorem}

\newtheorem{example}{Example}

\begin{document}
\vspace*{0.2in}

\begin{flushleft}
{\Large
\textbf\newline{
CRITERIA: A network decomposition and elementary flux mode translation-based tool for computing equilibria of biochemical systems
} 
}
\newline
\\
Exequiel Jun V. Villejo, 
Aurelio A. de los Reyes V, 
Bryan S. Hernandez\textsuperscript{*}
\\
\bigskip
Institute of Mathematics, University of the Philippines Diliman, Quezon City 1101, Philippines
\\
\bigskip

%
%





* bshernandez@up.edu.ph


\end{flushleft}
\section*{Abstract}
Analytically deriving equilibria, sometimes referred to as steady states, which often govern the long-term behavior of biochemical reaction networks, is essential for understanding cellular decision-making and robustness, yet remains computationally challenging for large and complex systems. The COMPILES framework of Hernandez et al. addresses this problem via network decomposition but is limited by a computationally intensive translation step and by its restriction to networks with zero kinetic deficiency. We introduce CRITERIA (Computing paRametrized posITive EquilibRIA), a new framework that overcomes these limitations through two key advances. First, it replaces the translation mechanism of COMPILES with a more efficient graph-theoretic formulation based on the work of Johnston and Burton, which generates a reaction-to-reaction graph from elementary flux modes and identifies directed cycles of the constructed graph via a binary linear program. Second, CRITERIA computes equilibria on a single unified translated network rather than solving subnetworks independently, thereby eliminating interdependencies that previously required extensive symbolic manipulation. Across a benchmark set of 26 biochemical models, CRITERIA achieves consistent and often substantial speed improvements while expanding applicability beyond the restricted class of zero kinetic deficiency systems. We demonstrate the biological utility of the framework by analyzing multistationarity, which underlies cellular decision-making, and absolute concentration robustness, a key mechanism for maintaining stable biochemical outputs, in the EnvZ–OmpR signaling pathway and a large-scale CRISPRi toggle switch, respectively. By improving both scalability and applicability, CRITERIA enables a systematic equilibrium analysis in biochemical networks of realistic size and complexity, providing a practical tool for studying long-term dynamical behavior in systems biology.

\section*{Author summary}
Understanding how biochemical systems settle into stable states, such as how protein concentrations reach equilibrium, is central to explaining cellular behavior and designing synthetic biological circuits. However, existing analytical tools for computing these equilibria, such as COMPILES, are limited by computational bottlenecks and can only be applied to a restricted class of reaction networks.
In this work, we introduce CRITERIA (Computing paRametrized posITive EquilibRIA), a new computational framework that makes equilibrium analysis more efficient and broadly applicable. CRITERIA uses a graph-based approach built on elementary flux modes to streamline key steps in the computation. It also changes how the problem is solved by combining subnetworks into a single system before computing equilibria, which avoids complicated symbolic calculations required in previous methods.
We demonstrate the usefulness of CRITERIA by studying biologically important systems, including the EnvZ-OmpR signaling pathway and a synthetic CRISPRi circuit. Our approach enables faster and more scalable analysis, allowing researchers to better understand how complex biochemical networks behave over time.


\section*{Introduction}
An equilibrium (or steady state) of a system refers to a state in which the concentrations of species in the system remain constant over time \cite{Allen2007,HornJackson1972}. 
The study of equilibria is fundamental to understanding biochemical systems, which are typically modeled using ordinary differential equations based on mass-action kinetics, since they characterize the long-term behavior of these systems \cite{Price2007,HernandezPLOS2023}. They reveal
important dynamical properties of the system, such as multistationarity, which refers to the capacity of the system to admit multiple positive equilibria, and absolute concentration robustness, where the concentration of a particular species remains the same at every positive equilibrium regardless of initial conditions.

Recent research has focused on network-based approaches in derivering analytic positive equilibria of biochemical systems \cite{Millan2012,Johnston2014,Thomson2009,JMP2019:parametrization,HernandezPLOS2023}. In particular, a publicly-available computational package built in MATLAB called COMPILES (COMPutIng anaLytic stEady States) was developed to facilitate this process \cite{HernandezPLOS2023}. Chemical reaction network (CRN) theory is the main analytical framework used by COMPILES, wherein the structural properties of \textit{weak reversibility} (WR) and \textit{deficiency zero} (DZ) are exploited to derive analytic equilibria. A network is WR if each of its reactions belongs to a directed cycle, and is said to be DZ if its deficiency, which is a nonnegative integer that measures the linear dependency of its reactions, is equal to zero \cite{FeinbergBook2019}. 

The method in COMPILES begins with \textit{network decomposition}, a divide and conquer approach to efficiently handle large networks, by decomposing the network into its \textit{independent subnetworks} \cite{FeinbergBook2019, HernandezPLOS2023, HDLC2021}. This is a special kind of network decomposition with the property that the rank of the parent network is equal to the sum of the ranks of its subnetworks. This is then  followed by \textit{network translation}, a process through which each subnetwork that initially does not possess WR and DZ is transformed into one that satisfies these properties while preserving the original dynamics of the system \cite{Johnston2014,HongSIAM}. The translation scheme implemented in COMPILES follows the method in \cite{HongSIAM}, in which pre- and post-translational filters are applied to identify, from a large number of candidates, translated networks with WR and DZ. Network translation produces a \textit{generalized chemical reaction network} (GCRN), which consists of two associated structures: \textit{stoichiometric CRN} and \textit{kinetic-order CRN}. The stoichiometric CRN has vertices corresponding to the translated complexes (i.e., reaction nodes), whereas the kinetic-order CRN retains vertices corresponding to the original complexes of the original CRN. Once the stoichiometric and kinetic-order CRNs have been successfully translated, equilibria are then parametrized per subnetwork based on \cite{JMP2019:parametrization}. Finally, the equilibria of each subnetwork are merged to obtain the equilibria of the entire network. This is guaranteed because independent decomposition has the property that the intersection of the set of positive equilibria of the subnetworks is equal to that of the entire network. \cite{FeinbergBook2019, HDLC2021}. In \cite{HernandezPLOS2023}, a sufficient condition was established to ensure that this property is maintained for any choice of rate constants. 

Despite its potential, COMPILES has several issues that make the analysis of biochemical reaction networks via equilibrium parametrization challenging. First, the network translation scheme employed by COMPILES constitutes a major bottleneck in the overall procedure. Although filters are incorporated to construct networks that satisfy WR and DZ, the exhaustive strategy of generating all possible candidate network translates imposes a significant computational overhead, making COMPILES impractical for larger systems. Second, there is a merging issue in COMPILES resulting to interdependencies in the derived equilibrium parametrization. In other words, this merging issue manifests in the presence of multiple species within each other’s equilibrium expressions. For instance, consider the following equilibrium parametrization of two species in the EnvZ–OmpR model \cite{JMP2019:parametrization}
\begin{align*}
    XD = \dfrac{X_pY k_8 k_{12}}{k_5\sigma_1} \hspace{1cm} X_pY = \dfrac{XDk_1k_3k_5}{k_2k_8(k_4+k_5)}
\end{align*}
where $k_i$'s are rate constants and $\sigma_1$ is a free parameter. It can be observed that both species appear in each other's equilibrium expressions. Merging equilibrium expressions requires extensive symbolic manipulation, which may cause MATLAB to fail in obtaining a closed-form solution free of interdependencies. Although not intrinsic to COMPILES, it nonetheless limits the method's applicability to more complex biochemical networks. Lastly, COMPILES only solves the equilibria of reaction networks whose kinetic deficiency (i.e., the deficiency of the kinetic-order CRN) is zero.

In this paper, we present an enhanced framework for computing the equilibria (i.e., deriving the analytic steady states) of biochemical systems by addressing the limitations of COMPILES as discussed above. We develop a MATLAB-based computational package called CRITERIA (Computing paRametrized posITive EquilibRIA) to facilitate this improved process. Here, we replace the translation scheme of COMPILES with the method developed in \cite{JohnstonBurton2019} to construct WR and DZ translations using an elementary flux mode-based approach. We also modify the order of steps in COMPILES to resolve the merging issue. In particular, rather than solving for the positive equilibria of each subnetwork independently, we first merge the translated (and, if there are any, untranslated) subnetworks into a single network and then derive the equilibrium parametrization of the whole network. Refer to Fig \ref{fig:generalmethod} for a more detailed discussion on the overall procedure of CRITERIA illustrated using a simple example. Finally, we incorporate the additional conditions required to parametrize the equilibria of networks with positive kinetic deficiency as established in \cite{JMP2019:parametrization}.

To demonstrate how CRITERIA improves COMPILES, we compare the runtime performance of both computational frameworks across a wide range of models. We further demonstrate the utility of the enhanced procedure by examining critical dynamical properties of complex biochemical systems, including absolute concentration robustness (ACR) in the full CRISPRi toggle switch model and the capacity for multistationarity in the  EnvZ-OmpR signaling pathway.



\FloatBarrier
\section*{Results}

\subsection*{An elementary flux mode-based approach for constructing weakly reversible and deficiency zero network translation}
\label{sec: translation}

One of the improvements made to resolve a key issue of COMPILES is the replacement of its network translation scheme with a more efficient procedure. In \cite{JohnstonBurton2019}, a computational method was developed to construct a WR and DZ translation via an elementary flux mode-based approach. \textit{Elementary flux modes} (EFMs) are flux-balanced pathways in the network that cannot be simplified in the sense that it is not possible to remove a subset of active reactions from an EFM and still be able to build a flux-balanced path using only the remaining active reactions \cite{JohnstonBurton2019,KlapperSIAM2021} (see the Methods section for more details).

The approach proposed in \cite{JohnstonBurton2019} consists of three main steps:
\begin{enumerate}
    \item From the given CRN, compute the set of EFMs $\mathcal{E} = \{e_1,\ldots,e_p\}$.
    \item Construct a reaction-to-reaction graph $G^\mathcal{R} = \left(V^\mathcal{R},E^\mathcal{R}\right)$ which is common source compatible and elementary flux mode compatible with the given CRN.
    \item Determine the translation complexes to produce the WR and DZ network translation.
\end{enumerate}

To illustrate the process, we consider the following histidine kinase system below \cite{Conradi2017}. Here, the histidine kinase ($X$) can autophosphorylate ($X_p$) and can transfer the phosphate group to a response regulator ($Y$) yielding ($Y_p$), which undergoes autodephosphorylation.

\begin{center}
\begin{tikzpicture}[baseline=(current  bounding  box.center)]
\tikzset{vertex/.style = {draw=none,fill=none}}
\tikzset{edge/.style = {->,> = latex', line width=0.10mm}} 
\node[vertex] (1) at  (-5,0) {$X$};
\node[vertex] (2) at  (-3,0) {$X_p$};
\node[vertex] (3) at  (-1.5,0) {$X_p +Y$};
\node[vertex] (4) at  (1.5,0) {$X + Y_p$};
\node[vertex] (5) at  (3,0) {$Y_p$};
\node[vertex] (6) at  (5,0) {$Y$};
\draw [edge]  (1) to["$k_1$"] (2);
\draw [edge, bend left = 10]  (3) to["$k_2$"] (4);
\draw [edge, bend left = 10]  (4) to["$k_3$"] (3);
\draw [edge]  (5) to["$k_4$"] (6);
\end{tikzpicture}
\end{center}

\textbf{Step 1.} Since EFMs represent flux-balanced pathways in the network, it is necessary to track the net production and consumption of species in the CRN. This information is codified in the \textit{stoichiometric matrix} $N$ (see the Methods section for more details on the basics of chemical reaction networks). Given this, an EFM is a flux vector that satisfies $Nv=0$ together with the nonnegative flux constraint $v \geq 0$. The central mathematical constructs that describe the space of all flux vectors that satisfy $Nv = 0$ where $v \geq 0$ (i.e., admissible flux vectors) are polyhedral cones \cite{KlapperSIAM2021}.

By defining $A = [I_r \hspace{.2cm} N \hspace{.1cm} -N]^T$, we can observe that the set
$$P = \{v \in \mathbb{R}^r_{\geq0} \hspace{.1cm} | \hspace{.1cm} Av \geq 0 \}$$
is a polyhedral cone that contains exactly the admissible flux vectors \cite{KlapperSIAM2021}. Moreover, it has been established in \cite{KlapperSIAM2021} that the extreme rays of this polyhedral cone correspond to the EFMs of the network. Alternatively, the set of EFMs is the minimal generating set of $P$.

The general method that computes this minimal generating set or the extreme rays of $P$ is called the \textit{Double Description Algorithm}. Several variations of the algorithm exists. In this paper, we use the \textit{Canonical Basis Method} which exploits the structure of the matrix $A$ as defined above \cite{KlapperSIAM2021}. The outline of the method is presented in Algorithm \ref{algo:CanonicalBasis}.

    
    
    
    


\begin{algorithm}
\caption{Computing elementary flux modes}
\label{algo:CanonicalBasis}
\begin{algorithmic}
\State \textbf{Input:} Stoichiometric matrix $N$
\State \textbf{Output:} A matrix $R$ whose columns are the elementary flux modes

\State \textbf{1. Initialize:} Set $k=r+2$, $A_\text{old} = I_r$, and $R_\text{old} = I_r$.

\State \textbf{2. Prepare rays in $R_\text{old}$ for partitioning:} While $k < r+2m+2$, do:
\Statex \hspace{1.5em} a. Select a row $v$ of $N$ not in $A_\text{old}$.
\Statex \hspace{1.5em} b. Update $A_\text{new} = \left[ A_\text{old} \hspace{.1cm} ; \hspace{.1cm} v \hspace{.1cm} ; \hspace{.1cm} -v \right]$.
\Statex \hspace{1.5em} c. Let $R_\text{new} = R_\text{old}$, and $\Lambda = \{1,2,\ldots,q\}$ where $q$ is the number of columns of $R_{\text{old}}$.

\State \textbf{3. Partition rays in $R_\text{old}$:} Let $\Lambda_+ = \{ j \in \Lambda : v \cdot R_\text{old}(:,j) > 0 \}$, $\Lambda_0 = \{ j \in \Lambda : v \cdot R_\text{old}(:,j) = 0 \}$, and $\Lambda_- = \{ j \in \Lambda : v \cdot R_\text{old}(:,j) < 0 \}$.

\State \textbf{4. Generate extreme rays:} For each $(j_+,j_-) \in \Lambda_+ \times \Lambda_-$, do:
\Statex \hspace{1.5em} a. Compute $r_{j_+,j_-} = (v \cdot R_{\text{old}}(:, j_+)) R_{\text{old}}(:, j_-) - (v \cdot R_{\text{old}}(:, j_-)) R_{\text{old}}(:, j_+)$.
\Statex \hspace{1.5em} b. Test for extremeness. If $r_{j_+,j_-}$ is an extreme ray, append $r_{j_+,j_-}$ to $R_\text{new}$.

\State \textbf{5. Update $R_\text{old}$ and $R_\text{new}$:}
\Statex \hspace{1.5em} a. Remove columns of $R_\text{new}$ indexed by $\Lambda_+ \times \Lambda_-$.
\Statex \hspace{1.5em} b. Update: $R_\text{old} = R_\text{new}$, $A_\text{old} = A_\text{new}$, and $k=k+2$.

\State \textbf{6. Repeat:} Go back to Step 2. Stop when no row of $N$ remains unchosen.
\end{algorithmic}
\end{algorithm}

Going back to the histidine kinase system consisting of four reactions, we can define its stoichiometric matrix $N$ as follows
\begin{equation*}
    \begin{blockarray}{ccccc}
     & R_1 & R_2 & R_3 & R_4 \\
    \begin{block}{c[rrrr]}
    X & -1 & 1 & -1 & 0 \\ 
    X_p & 1 & -1 & 1 & 0 \\
    Y & 0 & -1 & 1 & 1 \\
    Y_p & 0 & 1 & -1 & -1 \\
\end{block}
\end{blockarray}.
\end{equation*}

Applying Algorithm \ref{algo:CanonicalBasis}, the set of EFMs for the network are $e_1 = \{R_2,R_3\}$ and $e_2 = \{R_1, R_2, R_4\}$. See the Supplementary Information for the step-by-step procedure in the computation.

\textbf{Step 2.} Once the EFMs are computed, a reaction-to-reaction graph is then constructed that must be both common source (CS) compatible and elementary flux mode (EM) compatible with the given CRN (see Theorem \ref{thm: structuraltranslation} in the Methods section). 

A \textit{reaction-to-reaction graph} is a directed graph where the reactions of the CRN are treated as the vertices of the graph \cite{JohnstonBurton2019}. The CRN is said to be \textit{CS-compatible} whenever reactions of the following form in the CRN

\begin{center}
\begin{tikzpicture}[baseline=(current  bounding  box.center)]
\tikzset{vertex/.style = {draw=none,fill=none}}
\tikzset{edge/.style = {->,> = latex', line width=0.10mm}} 
\node[vertex] (3) at  (-1,0) {$\cdots$};
\node[vertex] (4) at  (.5,0) {$y$};
\node[vertex] (5) at  (1.5,.5) {};
\node[vertex] (6) at  (1.5,-.5) {};
\draw [edge]  (3) to["$r_k$"] (4);
\draw [edge]  (4) to["$r_i$"] (5);
\draw [edge]  (4) to node[below] {$r_j$} (6);
\end{tikzpicture}
\end{center}
forces $(r_k,r_i) \in E^\mathcal{R}$ if and only if $(r_k,r_j) \in E^\mathcal{R}$, where $E^\mathcal{R}$ is the edge set of the reaction-to-reaction graph \cite{JohnstonBurton2019}. Moreover, both the CRN and the reaction-to-reaction graph are said to be \textit{EM-compatible} if every EFM of the CRN corresponds to the vertices of a minimal directed cycle in the reaction-to-reaction graph, and vice versa \cite{JohnstonBurton2019}.

In \cite{JohnstonBurton2019}, a binary linear programming (BLP) problem was formulated to construct the desired reaction-to-reaction graph while enforcing CS compatibility and EM compatibility with the given CRN. In our histidine kinase system, the BLP model results to the reaction-to-reaction graph shown below. See the Supplementary Information for the detailed discussion.

\begin{center}
\begin{tikzpicture}[baseline=(current  bounding  box.center)]
\tikzset{vertex/.style = {draw=none,fill=none}}
\tikzset{edge/.style = {->,> = latex', line width=0.10mm}} 
\node[vertex] (1) at  (1.5,-1.5) {$R_1$};
\node[vertex] (2) at  (0,0) {$R_2$};
\node[vertex] (3) at  (-1.5,-1.5) {$R_3$};
\node[vertex] (4) at  (0,-2.5) {$R_4$};
\draw [edge]  (2) to (1);
\draw [edge]  (4) to (2);
\draw [edge]  (1) to (4);
\draw [edge, bend left = 10]  (2) to (3);
\draw [edge, bend left = 10]  (3) to (2);
\end{tikzpicture}
\end{center}

\textbf{Step 3.} The final step of the translation procedure is to determine the translation complexes to obtain a WR and DZ network translation. Intuitively, \textit{translation complexes} are complexes added to individual reactions in the translation process. This modifies the graphical structure of the CRN without changing its stoichiometric matrix \cite{Johnston2014}. This is crucial for preserving the dynamics of the system.

Theorem \ref{thm: structuraltranslation} provides the method for determining the translation complexes, which involves solving the linear system (\ref{eq:2}). A more efficient procedure, however, is proposed in \cite{JohnstonBurton2019} to determine the translation complexes instead of solving the linear system (\ref{eq:2}) directly. 

The approach is similar to a graph traversal: once a value is fixed at a single node, the remaining values are computed incrementally using Eq (\ref{eq:2}) by traversing the directed edges of the reaction-to-reaction graph. The outline of the procedure is shown in Algorithm \ref{algo:TranslationComplexes}. This strategy exploits the structure of the reaction-to-reaction graph to propagate values of the translation complexes step-by-step. This is more efficient than solving the full matrix system which may have redundancies that can lead to unnecessary computational overhead.

\begin{algorithm}
\caption{Determining the translation complexes}
\label{algo:TranslationComplexes}
\begin{algorithmic}
\State \textbf{Input:} Edges $E^\mathcal{R}$ of the reaction-to-reaction graph, matrix of source complexes $Y_s$ and product complexes $Y_p$, linkage classes $\mathcal{L}$
\State \textbf{Output:} Translation complexes $\alpha$ \\
    \State \textbf{1. Check consistency of the linear system.} Determine whether $\text{rank}(A) = \text{rank}(B)$, where $B$ is the augmented matrix $B = [A \hspace{.1cm} | \hspace{.1cm} b]$, and $A$ and $b$ come from Eq (\ref{eq:2}). Exit algorithm if $\text{rank}(A) \neq \text{rank}(B)$.
    
    \State \textbf{2. Seed the process.} Set $\mathcal{P} = \mathcal{R}$. Select an arbitrary $r_i \in \mathcal{P}$ and set $\alpha_i = 0$. Move $r_i$ to the "active" list $\mathcal{P}'$.
    
    \State \textbf{3. Expand.} For every $r_i$ in the "active" list, choose a neighbor $r_j \in \mathcal{P}$, and then calculate $\alpha_j$ based on Eq (\ref{eq:2}). Move $r_j$ to the "processed" list $\mathcal{P}''$.
    
    \State \textbf{4. Propagate.} Once we have checked all neighbors of the current "active" list $\mathcal{P}'$, set $\mathcal{P}' = \mathcal{P}''$. Set $\mathcal{P}'' = \emptyset$ and repeat Step 3.
    
    \State \textbf{5. Check for completion.} If the "processed" list is already empty but $\mathcal{P} \neq \emptyset$, go back to step 2. If $\mathcal{P} = \emptyset$, exit the algorithm.

    \State \textbf{6. Update complexes.} Add the computed translation complexes $\alpha_i$ to reaction $i$.

    \State \textbf{7. Ensure nonnegative complexes.} For each linkage class, shift complexes with negative entries to ensure that the stoichiometric coefficients of species per complex is nonnegative. 
\end{algorithmic}
\end{algorithm}

In our running example, the histidine kinase network consists of four reactions, and so four translation complexes must be determined. Applying Algorithm \ref{algo:TranslationComplexes}, the translation complex for each reaction is given by
\begin{align*}
    X &\longrightarrow X_p \hspace{1.5cm} (+Y_p) \\
    X_p + Y &\longrightarrow X+Y_p \hspace{.8cm} (+\emptyset) \\
    X + Y_p &\longrightarrow X_p+Y \hspace{.75cm} (+\emptyset) \\
    Y_p &\longrightarrow Y \hspace{1.63cm} (+X_p) \\
\end{align*}
which results in the following translated network which is guaranteed to be WR and DZ by Theorem \ref{thm: structuraltranslation}. See the Supplementary Information for the step-by-step procedure in the computation.

\begin{center}
    \begin{tikzpicture}[baseline=(current  bounding  box.center)]
    \tikzset{vertex/.style = {draw=none,fill=none}}
    \tikzset{edge/.style = {->,> = latex', line width=0.10mm}} 
    \node[vertex] (1) at  (-2,0) {$X+Y_p$};
    \node[vertex] (2) at  (2,0) {$X_p+Y_p$};
    \node[vertex] (3) at  (0,-1.5) {$X_p+Y$};
    \draw [edge]  (1) to["$k_1$"] (2);
    \draw [edge, bend left = 10]  (3) to["$k_2$"] (1);
    \draw [edge, bend left = 10]  (1) to["$k_3$"] (3);
    \draw [edge]  (2) to["$k_4$"] (3);
    \end{tikzpicture}
    \end{center}

It is worth noting that a WR and DZ translation from the method of Theorem \ref{thm: structuraltranslation} is not unique. For instance, translating using $\alpha_1 = Y$, $\alpha_2=\emptyset=\alpha_3$, and $\alpha_4 = X$ yields a different network but is nevertheless WR and DZ. This arises from Steps 2 and 3 of Algorithm \ref{algo:TranslationComplexes}, where different choices of $r_i$ from $\mathcal{P}$ may lead to different outcomes.

In the next section, we demonstrate the enhanced computational framework for deriving the analytic steady states using a two-protein gene transcription model.

\subsection*{Demonstration: An enhanced framework for deriving the analytic steady state of a two-protein gene transcription model}

We now demonstrate the enhanced framework for deriving analytic equilibria by considering a two-protein transcription model, which is a mass-action system consisting of ten reactions that describe a gene transcription motif involving two proteins $P_1$ and $P_2$, each produced by their respective genes $G_1$ and $G_2$ \cite{Conradi2017,SiegalGaskins2015}. Table \ref{tab:2protein} shows the reactions of the network grouped according to their functional roles.

\begin{center}
\renewcommand{\arraystretch}{1.5}
\begin{longtable}{lll}

\toprule
Label & Reaction & Description \\
\midrule
\endfirsthead

\toprule
Label & Reaction & Description \\
\midrule
\endhead



\multicolumn{3}{l}{\textbf{Gene transcription}} \\
$R_1$ & $G_1 \xrightarrow{k_1} G_1 + P_1$ & basal production of protein $P_1$ through gene $G_1$ \\
$R_2$ & $G_2 \xrightarrow{k_2} G_2 + P_2$ & basal production of protein $P_2$ through gene $G_2$ \\

\midrule
\multicolumn{3}{l}{\textbf{Protein clearance}} \\
$R_3$ & $P_1 \xrightarrow{k_3} 0$ & clearance of protein $P_1$ \\
$R_4$ & $P_2 \xrightarrow{k_4} 0$ & clearance of protein $P_2$ \\

\midrule
\multicolumn{3}{l}{\textbf{Dimerization}} \\
$R_5$ and $R_6$ & $2P_2 \xrightleftharpoons[k_6]{k_5} D$ & dimerization of protein $P_2$ \\

\midrule
\multicolumn{3}{l}{\textbf{Regulatory complex formation}} \\
$R_7$ and $R_8$ & $G_2 + P_1 \xrightleftharpoons[k_8]{k_7} C_1$ & $P_1$ binds with $G_2$ \\
$R_9$ and $R_{10}$ & $G_1 + D \xrightleftharpoons[k_{10}]{k_9} C_2$ & $P_2$ dimer binds with $G_1$ \\ 
\bottomrule
\caption{\scriptsize {\bf The reactions of the two-protein gene transcription model and their descriptions.}}
\label{tab:2protein}
\end{longtable}
\end{center}

\vspace{-.5cm}

\begin{figure*}
\centering
\includegraphics[width=13.5cm,height=75.5cm,keepaspectratio]{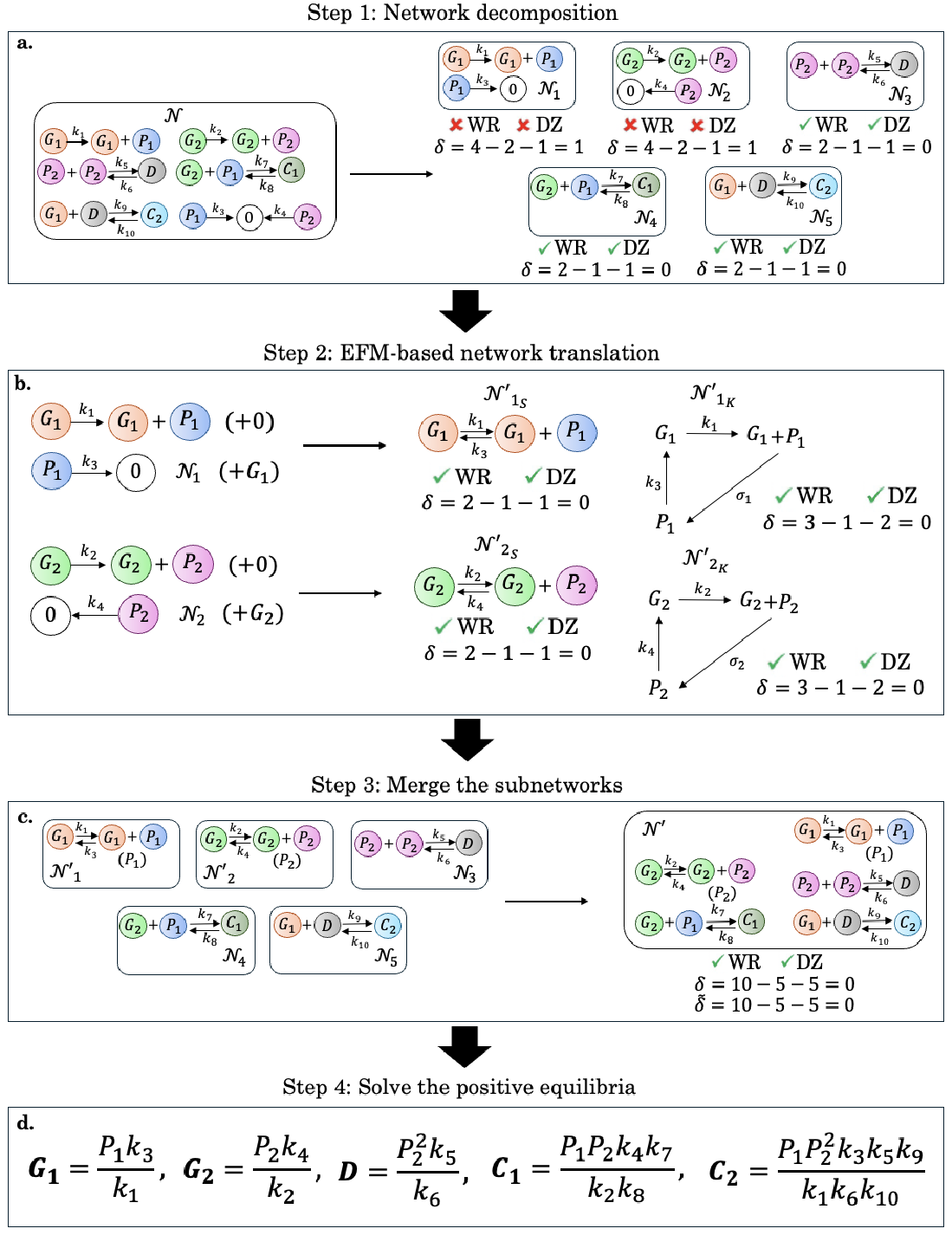}
\caption{\scriptsize
{\bf{Derivation of the analytic equilibria using the enhanced framework.}}
{\bf{a}} The CRN $\mathcal{N}$ composed of $10$ reactions is decomposed into five independent subnetworks, each having two reactions. This decomposition has the property that the rank of the stoichiometric matrix of $\mathcal{N}$ is equal to the sum of the ranks of the stoichiometric matrices of the five subnetworks.
{\bf{b}} Of the five subnetworks, two ($\mathcal{N}_1$ and $\mathcal{N}_2$) are translated to become WR and DZ. For $\mathcal{N}_1$, the complex $G_1$ serves as the translation complex for the third reaction while the first reaction remains unchanged. This produces a stoichiometric CRN ($\mathcal{N}'_{1_S}$) and kinetic-order CRN ($\mathcal{N}'_{1_K}$) that are both WR and DZ. Hence, the translated network $\mathcal{N}'_1$ is WR and DZ, while retaining the same dynamics as $\mathcal{N}_1$. The same procedure is done for $\mathcal{N}_2$ by using $G_2$ as the translation complex for the fourth reaction and keeping the second reaction unchanged. This yields a translated network with the same dynamics as the original network but is now WR and DZ.
{\bf{c}} The five subnetworks are then merged into a single network. For $\mathcal{N}'_1$ and $\mathcal{N}'_2$, the complex in parentheses is the kinetic complex associated in that node. For nodes where no complex is shown in parentheses, it is assumed that the stoichiometric complex is the same as the kinetic complex. The kinetic complexes comprise the kinetic-order CRN, with additional edges, known as phantom edges, introduced to separate multiple complexes associated with a single node. There are two deficiencies associated to the merged network. The first one ($\delta$) is for the deficiency of the stoichiometric CRN while the second one ($\tilde{\delta}$) is for the kinetic-order CRN. In the example, both deficiencies are equal to zero and hence the network is DZ.
{\bf{d}} Since the network is now WR and DZ, the equilibrium expressions are then analytically derived.
}\label{fig:generalmethod}
\end{figure*}

The first step of the procedure is to break down the network into its independent subnetworks as shown in Fig \ref{fig:generalmethod}a. The rank of the entire network is given by
\begin{equation*}
    s = \text{rank} \begin{blockarray}{ccccccccccc}
     & R_1 & R_2 & R_3 & R_4 & R_5 & R_6 & R_7 & R_8 & R_9 & R_{10} \\
    \begin{block}{c[rrrrrrrrrr]}
    G_1 & 0 & 0 & 0 & 0 & 0 & 0 & 0 & 0 & -1 & 1 \\ 
    G_2 & 0 & 0 & 0 & 0 & 0 & 0 & -1 & 1 & 0 & 0 \\ 
    P_1 & 1 & 0 & -1 & 0 & 0 & 0 & -1 & 1 & 0 & 0 \\ 
    P_2 & 0 & 1 & 0 & -1 & -2 & 2 & 0 & 0 & 0 & 0 \\ 
    D & 0 & 0 & 0 & 0 & 1 & -1 & 0 & 0 & -1 & 1 \\ 
    C_1 & 0 & 0 & 0 & 0 & 0 & 0 & 1 & -1 & 0 & 0 \\ 
    C_2 & 0 & 0 & 0 & 0 & 0 & 0 & 0 & 0 & 1 & -1 \\ 
\end{block}
\end{blockarray} = 5
\end{equation*}
while the ranks of the subnetworks are as follows
\begin{equation*}
    s_1 = \text{rank} \begin{blockarray}{cc}
     R_1 & R_3 \\
    \begin{block}{[rr]}
    0 & 0 \\ 
    0 & 0 \\ 
    1 & -1 \\  
    0 & 0 \\  
    0 & 0 \\  
    0 & 0 \\  
    0 & 0 \\  
\end{block}
\end{blockarray} = 1 \hspace{.5cm} s_2 = \text{rank} \begin{blockarray}{cc}
     R_2 & R_4 \\
    \begin{block}{[rr]}
    0 & 0 \\ 
    0 & 0 \\ 
    0 & 0 \\  
    1 & -1 \\  
    0 & 0 \\  
    0 & 0 \\  
    0 & 0 \\  
\end{block}
\end{blockarray} = 1 \hspace{.5cm} s_3 = \text{rank} \begin{blockarray}{cc}
     R_5 & R_6 \\
    \begin{block}{[rr]}
    0 & 0 \\ 
    0 & 0 \\ 
    0 & 0 \\  
    -2 & 2 \\  
    1 & -1 \\  
    0 & 0 \\  
    0 & 0 \\  
\end{block}
\end{blockarray} = 1 
\end{equation*}

\begin{equation*}
    s_4 = \text{rank} \begin{blockarray}{cc}
     R_7 & R_8 \\
    \begin{block}{[rr]}
    0 & 0 \\ 
    -1 & 1 \\ 
    -1 & 1 \\  
    0 & 0 \\  
    0 & 0 \\  
    1 & -1 \\  
    0 & 0 \\  
\end{block}
\end{blockarray} = 1 \hspace{.5cm} s_5 = \text{rank} \begin{blockarray}{cc}
     R_9 & R_{10} \\
    \begin{block}{[rr]}
    -1 & 1 \\ 
    0 & 0 \\ 
    0 & 0 \\  
    0 & 0 \\  
    -1 & 1 \\  
    0 & 0 \\  
    1 & -1 \\  
\end{block}
\end{blockarray} = 1
\end{equation*}

Since $s = 5 = s_1+s_2+s_3+s_4+s_5$, the decomposition in Fig \ref{fig:generalmethod}a is an independent network decomposition. Among the five subnetworks, three are already WR and DZ, whereas the first two require network translation to satisfy these conditions. Following the translation scheme outlined in the previous section, we use $G_1$ as the translation complex for the third reaction, while leaving the first reaction unchanged. Refer to Fig \ref{fig:generalmethod}b. As a result, the third reaction becomes $G_1+P_1 \longrightarrow G_1$ which is the reverse of the first reaction. The translated subnetwork with reversible reactions defines the stoichiometric CRN, denoted by $\mathcal{N}'_{1_S}$, associated with $\mathcal{N}_1$. Clearly, it is both WR and DZ $(\delta' = n'-l'-s' = 2-1-1=0)$. 

Next, we construct the kinetic-order CRN, which is denoted by $\mathcal{N}'_{1_K}$, associated with $\mathcal{N}_1$. To do this, we first note that we keep the reaction rates of the original network so that the dynamics of the system do not change. In particular, we fix $k_3P_1$ as the reaction rate for $G_1 + P_1 \longrightarrow G_1$. Now, we take all the edges of $\mathcal{N}'_{1_S}$ while assigning to each of the source node of the edges the corresponding source complexes from the original network. Since $P_1$ is the original source complex of the third reaction, two complexes become associated with a single node in $\mathcal{N}'_{1_K}$. We separate them using a \textit{phantom edge}, defined as an edge equipped with a free parameter ($\sigma_1$) serving as its "rate constant", in such a way that the resulting network is WR. For a phantom edge, we always get a zero stoichiometric vector since both the head and tail of the edge correspond to the same complex. This construction yields the kinetic-order CRN $\mathcal{N}'_{2_K}$, which is both WR and DZ $(\tilde{\delta} = \tilde{n} - \tilde{l} - \tilde{s} = 3-1-2=0)$. The stoichiometric CRN $\mathcal{N}'_{1_S}$ and kinetic-order CRN $\mathcal{N}'_{1_K}$ constitute the \textit{translated network} for $\mathcal{N}_1$. Since both $\mathcal{N}'_{1_S}$ and $\mathcal{N}'_{1_K}$ are WR and DZ, the translated network then is WR and DZ. 

Similarly, we apply the procedure to $\mathcal{N}_2$ using $G_2$ as the translation complex for the fourth reaction while keeping the second reaction unchanged. This yields a translated network that is both WR and DZ as shown in Fig \ref{fig:generalmethod}b.

The next step is to merge the subnetworks into a single unified network instead of treating each subnetwork as an individual piece. This step is crucial for avoiding dependencies in the equilibrium expressions, as described in the Introduction section. It circumvents the need for extensive symbolic manipulations in merging equilibrium expressions, which can otherwise overwhelm MATLAB.

Finally, we compute the equilibrium parametrization using the full network rather than treating each subnetwork individually. We refer the readers to \cite{JMP2019:parametrization,HernandezPLOS2023} for a detailed discussion on equilibrium parametrization.

To summarize the improvements introduced in CRITERIA, we make three modifications in the framework of COMPILES. First, we adopt the method developed in \cite{JohnstonBurton2019} to translate networks into becoming WR and DZ using elementary flux modes. Second, we reorder the steps of COMPILES. Instead of computing the equilibrium parametrization of each subnetwork individually, we first merge all subnetworks into a single network by taking the union of their reactions, and then compute the parametrization for the entire network. Finally, we extend the applicability of the procedure to handle cases in which the deficiency of the kinetic-order CRN is positive, based on the theory established in \cite{JMP2019:parametrization}. This is in contrast to COMPILES which only considers networks with zero kinetic deficiency.

In the next section, we compare the runtime performance of COMPILES and CRITERIA by applying them to a variety of models.

\subsection*{Runtime performance comparisons}

We apply CRITERIA to 26 models and evaluate its performance by comparing its runtime with that of COMPILES. We also assess the performance gains of the new network translation scheme by comparing the EFM-based translation with the translation method used in COMPILES. Computational experiments were conducted on a MacBook Air equipped with an Apple M2 processor and 8 GB of RAM. All runtime measurements were obtained under identical hardware and software conditions.

To contextualize performance results, Table \ref{tab:networkproperties} summarizes the network properties of the benchmark models used for the comparisons. The models are listed in order of increasing number of reactions. It is woth noting here that all benchmark models are neither WR nor DZ.

\begin{center}
\renewcommand{\arraystretch}{1.2}
\begin{longtable}{p{7cm}ccccccc}

\toprule
Model & $r$ & $m$ & $n$ & $s$ & $l$ & $\delta$ & $\tilde{\delta}$ \\
\midrule
\endfirsthead

\toprule
Model & $r$ & $m$ & $n$ & $s$ & $l$ & $\delta$ & $\tilde{\delta}$ \\
\midrule
\endhead

Model A: Histidine-Kinase \cite{Conradi2017} & 4 & 4 & 6 & 2 & 3 & 1 & 0 \\
Model B: Miao's Influenza Virus Model \cite{FluModels2019} & 5 & 3 & 7 & 3 & 2 & 2 & 0 \\
Model C: IDHKP-IDH \cite{Shinar2009} & 6 & 5 & 6 & 3 & 2 & 1 & 0 \\
Model D: 1-Site PD Network \cite{Villareal2024} & 6 & 6 & 6 & 3 & 2 & 1 & 0 \\
Model E: Hybrid Histidine-Kinase \cite{Conradi2017} & 6 & 6 & 10 & 4 & 4 & 2 & 0 \\
Model F: PTM system \cite{Feliu2013} & 9 & 8 & 9 & 5 & 3 & 1 & 0 \\
Model G: 2-Protein Gene Transcription \cite{Conradi2017,SiegalGaskins2015} & 10 & 7 & 13 & 5 & 6 & 2 & 0 \\
Model H: Reduced ERK Network \cite{Obatake2019} & 10 & 10 & 12 & 7 & 2 & 3 & 0 \\
Model I: 2-site PD Network \cite{Villareal2024} & 12 & 9 & 10 & 6 & 2 & 2 & 0 \\
Model J: Two-Substrate Phosphorylation \cite{Conradi2017} & 12 & 10 & 12 & 6 & 4 & 2 & 0 \\
Model K: Two-Layer Cascade of Modification Cycles \cite{Feliu2013} & 12 & 10 & 12 & 6 & 4 & 2 & 0 \\
Model L: EnvZ-OmpR \cite{JMP2019:parametrization} & 14 & 9 & 13 & 7 & 4 & 2 & 1 \\
Model M: Two-Component Osmoregulator in \textit{E. coli} \cite{Karp2012} & 15 & 9 & 13 & 7 & 4 & 2 & 2 \\
Model N: Hernandez's Influenza Virus Model \cite{FluModels2019}  & 16 & 7 & 16 & 7 & 3 & 6 & 4 \\
Model O: Modified Two-Component Osmoregulator in \textit{E. coli} \cite{Karp2012} & 16 & 9 & 15 & 7 & 5 & 3 & 2 \\
Model P: CRISPRi Model 1 \cite{SantosMoreno2020}  & 16 & 10 & 19 & 8 & 9 & 2 & 0 \\
Model Q: MAPK \cite{Johnston2019} & 16 & 11 & 12 & 8 & 2 & 2 & 1 \\
Model R: Full ERK model \cite{Obatake2019} & 18 & 12 & 14 & 9 & 2 & 3 & 1 \\
Model S: 3-Site PD Network \cite{Villareal2024} & 18 & 12 & 14 & 9 & 2 & 3 & 1 \\
Model T: Zigzag Model of Plant-Pathogen Interactions \cite{Johnston2019} & 21 & 13 & 23 & 9 & 6 & 8 & 0 \\
Model U: CRISPRi Model 2 \cite{SantosMoreno2020} & 22 & 13 & 24 & 11 & 11 & 2 & 0 \\
Model V: MAPK Signaling Cascade \cite{Feliu2013} & 24 & 17 & 19 & 12 & 4 & 3 & 0 \\
Model W: WNT Signaling Pathway \cite{JMP2019:parametrization} & 31 & 19 & 28 & 14 & 10 & 4 & 0 \\
Model X: Insulin Signaling Pathway \cite{LML2022} & 35 & 20 & 35 & 15 & 13 & 7 & 1 \\
Model Y: CRISPRi Model 3 \cite{SantosMoreno2020}  & 38 & 17 & 34 & 15 & 13 & 6 & 8 \\
Model Z: CRISPRi Model 4 \cite{SantosMoreno2020} & 42 & 17 & 38 & 15 & 13 & 10 & 10 \\
\bottomrule
\caption{\scriptsize {\bf Network properties of benchmark models with $r=$ number of reactions, $m=$ number of species, $n=$ number of complexes, $s=$ rank of the network, $l=$ number of linkage classes, $\delta=$ effective deficiency, and $\tilde{\delta}$ = kinetic deficiency}}
\label{tab:networkproperties}
\end{longtable}
\end{center}

\vspace{-.5cm}

We impose a maximum runtime of 3600 seconds; any execution that does not complete within this time limit is terminated. In our first computational experiment, we compare the runtime performance of the EFM-based network translation, which we call REWARDZ (tRanslation via Elementary flux modes
towards WeAkly Reversible Deficiency Zero networks), with TOWARDZ (TranslatiOn toward
WeAkly Reversible and Deficiency Zero networks), the translation method used in COMPILES \cite{HongSIAM}. Table \ref{tab:towardzVSrewardz} presents this comparison. It should be noted that network decomposition has not yet been integrated in both computational frameworks. We make the following observations:
\begin{itemize}
    \item REWARDZ successfully produced a valid translation for nearly all models, with the exception of Hernandez’s Influenza Virus Model and CRISPRi Model 3, whereas TOWARDZ succeeded in only five models. In particular, TOWARDZ failed to translate any model containing more than six reactions.
    \item Among the models for which REWARDZ successfully constructed a WR and DZ translation, the runtime was under 3 seconds in all cases except for the CRISPRi full model (CRISPRi Model 4). This model, which is the largest considered in the paper, contains 42 reactions and has deficiency equal to 10.
    \item There is no significant difference in runtime among the models for which both methods successfully produced a valid translation, with the exception of the Hybrid Histidine–Kinase Model. This model, which consists of six reactions and ten species and has a deficiency of two, exhibited markedly better performance under REWARDZ.
\end{itemize}

The result in Table \ref{tab:towardzVSrewardz} demonstrates that REWARDZ offers a substantial improvement in terms of running time in generating a WR and DZ network translation, especially in larger models. It underscores the effectiveness and scalability of the enhanced computational framework in handling more complex biochemical systems.

\begin{center}
\renewcommand{\arraystretch}{1.2}
\begin{longtable}{c|c c|c c}

\toprule
\multicolumn{1}{c}{Model} & 
\multicolumn{1}{c}{TOWARDZ} & 
\multicolumn{1}{c}{REWARDZ} & 
\multicolumn{1}{c}{COMPILES} & 
\multicolumn{1}{c}{CRITERIA} \\
\midrule
\endfirsthead

\toprule
\multicolumn{1}{c}{Model} & 
\multicolumn{1}{c}{TOWARDZ} & 
\multicolumn{1}{c}{REWARDZ} & 
\multicolumn{1}{c}{COMPILES} & 
\multicolumn{1}{c}{CRITERIA} \\
\midrule
\endhead

Model A & 0.2594 & 0.4982 & 0.7623  & 0.9914  \\
Model B & 1.0428 & 0.5408 & 0.8771  & 0.9932  \\
Model C & 0.5568 & 0.6829 & 1.1345  & 1.2003  \\
Model D & 1.0404 & 0.7311 & 1.5676 & 1.1916 \\
Model E & 204.2957 & 0.7492 & 206.5428 & 1.3896 \\
Model F & $\ggg$ 3600 & 1.0984 & $\ggg$ 3600 & 1.8361 \\
Model G & $\ggg$ 3600 & 0.8535 & 1.0387 & 1.4518 \\
Model H & $\ggg$ 3600 & 1.1990 & $\ggg$ 3600 & 2.0491 \\
Model I & $\ggg$ 3600 & 1.2163 & 2.9069 & 2.0675 \\
Model J & $\ggg$ 3600 & 1.2474 & 2.7384 & 1.8826 \\
Model K & $\ggg$ 3600 & 1.2224 & 2.8884 & 1.8970 \\
Model L & $\ggg$ 3600 & 1.2796 & $\ggg$ 3600 & 2.5600  \\
Model M & $\ggg$ 3600 & 1.3573 & $\ggg$ 3600 & 2.7157  \\
Model N & $\ggg$ 3600 & No valid translation & $\ggg$ 3600 & 3.0672 \\
Model O & $\ggg$ 3600 & 1.4710 & $\ggg$ 3600 & 3.1765 \\
Model P & $\ggg$ 3600 & 1.0709 & 1.4306 & 2.0085 \\
Model Q & $\ggg$ 3600 & 1.1579 & 1238.6329 & 2.6412 \\
Model R & $\ggg$ 3600 & 1.3966 & $\ggg$ 3600 & 3.1966 \\
Model S & $\ggg$ 3600 & 1.4828 & Error & 2.9613 \\
Model T & $\ggg$ 3600 & 1.4232 & 15.8447 & 3.0097 \\
Model U & $\ggg$ 3600 & 1.3345 & Error & 2.4981 \\
Model V & $\ggg$ 3600 & 1.9116 & $\ggg$ 3600 & 3.7726 \\
Model W & $\ggg$ 3600 & 2.2248 & $\ggg$ 3600 & 3.8027 \\
Model X & $\ggg$ 3600 & 2.1477 & 3.6752 & 5.0770 \\
Model Y & $\ggg$ 3600 & No valid translation & $\ggg$ 3600 & 6.4235 \\
Model Z & $\ggg$ 3600 & 32.4301 & $\ggg$ 3600 & 12.9456 \\
\bottomrule

\caption{\scriptsize {\bf Comparison of running times (in seconds) between TOWARDZ and REWARDZ, and between COMPILES and CRITERIA.} A maximum runtime of 3600 seconds is imposed for the computational experiments. Entries marked with “$\ggg$ 3600” indicate that the method was unable to produce the desired output within the time limit. Entries labeled “No valid translation” signify that the method completed within the prescribed maximum runtime but did not identify a WR and DZ translation. Lastly, entries labeled “Error” encountered an error in the code.} 
\label{tab:towardzVSrewardz}
\end{longtable}
\end{center} 

\vspace{-.8cm}

We now proceed to compare the runtime performance of COMPILES and CRITERIA. Similar to the comparison above, we impose a maximum runtime of 3600 seconds. If steady states are not parametrized within the time limit, we exit the program. As shown in Table \ref{tab:towardzVSrewardz}, COMPILES was unable to derive the analytic equilibria of several benchmark models. Of the 26 models, it failed to parametrize the equilibria of 11, whereas CRITERIA successfully obtained equilibrium parametrizations for all models. Notably, CRITERIA was able to solve the analytic equilibria of all models in less than 7 seconds except for the full CRISPRi model which required approximately 12 seconds. In addition, COMPILES encountered code errors in 2 models, namely the 3-site PD Network and CRISPRi Model 2. These results show the performance gains of CRITERIA over COMPILES. 

Now, it can be observed that COMPILES performed slightly better than CRITERIA in six models. This is because, in these models, network decomposition has broken down the network into WR and DZ subnetworks, eliminating the need for translation. To elucidate this, we consider the 2-Protein Gene Transcription Model, CRISPRi Model 1, and the Insulin Signaling Pathway. Fig \ref{fig:generalmethod} shows that the 2-Protein Gene Transcription Model consists of five independent subnetworks, three of which are already WR and DZ. The remaining two subnetworks that require translation each contain only two reactions. For the Insulin Signaling Pathway, five of the ten independent subnetworks are already WR and DZ. The remaining five subnetworks that require translation only have three reactions, four reactions, and two reactions in three subnetworks, respectively. Lastly, the CRISPRi Model 1 (smallest CRISPRi model considered in this paper) consists of eight independent subnetworks, none of which are WR or DZ, but each contains only two reactions.

Based on these results, COMPILES appears to handle smaller networks slightly more efficiently than CRITERIA. This is likely due to the fact that CRITERIA requires solving elementary flux modes and a binary linear program, which can be computationally more demanding. For more details, see the Discussion section.

We verify the correctness of the output produced by CRITERIA by substituting the derived equilibrium expressions into the model's ODE system. The time derivatives evaluate to zero, confirming that the solutions are indeed valid equilibria.

In the next section, we apply CRITERIA to investigate key dynamical properties of selected biochemical systems such as multistationarity and ACR.

\subsection*{Application: Using equilibrium parametrization to examine ACR and multistationarity}

We now apply our method to the EnvZ-OmpR system, a two-component signaling pathway consisting of the sensor kinase EnvZ, denoted by $X$, and the response regulator OmpR, denoted by $Y$ \cite{ShinarFeinberg2010}. In the network below that corresponds to Model L in Table \ref{tab:networkproperties}, both $X$ and $Y$ exist in phosphorylated forms, denoted by $X_p$ and $Y_p$, respectively. The sensor $X$ undergoes authophosphorylation through ATP (species $T$) binding and hydrolysis. The phosphorylated sensor $X_p$ catalyzes the transfer of the phosphoryl group to $Y$. Furthermore, $X$ dephosphorylates $Y_p$ using either ATP or ADP (species $D$) as a cofactor \cite{ShinarFeinberg2010}. In the model, the effects of both ATP and ADP are simultaneously considered. Experimental observations indicate that the system exhibits concentration robustness, which has been analytically verified by demonstrating that the equilibrium concentration of $Y_p$ depends solely on the rate constants. \cite{ShinarFeinberg2010,JMP2019:parametrization,Batchelor2003}.

The reaction network of the system consisting of $14$ reactions is shown below. \begin{center}
\begin{tikzpicture}[baseline=(current  bounding  box.center)]
\tikzset{vertex/.style = {draw=none,fill=none}}
\tikzset{edge/.style = {->,> = latex', line width=0.10mm}} 
\node[vertex] (1) at  (-4,0) {$XD$};
\node[vertex] (2) at  (-2,0) {$X$};
\node[vertex] (3) at  (0,0) {$XT$};
\node[vertex] (4) at  (2,0) {$X_p$};

\node[vertex] (5) at  (-4,-1.5) {$X_p + Y$};
\node[vertex] (6) at  (-1.5,-1.5) {$X_p Y$};
\node[vertex] (7) at  (0.5,-1.5) {$X + Y_p$};

\draw [edge] ([yshift=2pt]1.east) -- ([yshift=2pt]2.west) node[midway, above] {$k_1$};
\draw [edge] ([yshift=-2pt]2.west) -- ([yshift=-2pt]1.east) node[midway, below] {$k_2$};

\draw [edge] ([yshift=2pt]2.east) -- ([yshift=2pt]3.west) node[midway, above] {$k_3$};
\draw [edge] ([yshift=-2pt]3.west) -- ([yshift=-2pt]2.east) node[midway, below] {$k_4$};
\draw [edge]  (3) to node[above] {$k_5$} (4);

\draw [edge] ([yshift=2pt]5.east) -- ([yshift=2pt]6.west) node[midway, above] {$k_6$};
\draw [edge] ([yshift=-2pt]6.west) -- ([yshift=-2pt]5.east) node[midway, below] {$k_7$};
\draw [edge]  (6) to node[above] {$k_8$} (7);
\end{tikzpicture}
\end{center}

\begin{center}
\begin{tikzpicture}[baseline=(current  bounding  box.center)]
\tikzset{vertex/.style = {draw=none,fill=none}}
\tikzset{edge/.style = {->,> = latex', line width=0.10mm}} 
\node[vertex] (8) at  (-4,-3) {$XD + Y_p$};
\node[vertex] (9) at  (-1.5,-3) {$XDY_p$};
\node[vertex] (10) at  (1,-3) {$XD + Y$};

\node[vertex] (11) at  (-4,-4.5) {$XT + Y_p$};
\node[vertex] (12) at  (-1.5,-4.5) {$XTY_p$};
\node[vertex] (13) at  (1,-4.5) {$XT + Y$};

\draw [edge] ([yshift=2pt]8.east) -- ([yshift=2pt]9.west) node[midway, above] {$k_9$};
\draw [edge] ([yshift=-2pt]9.west) -- ([yshift=-2pt]8.east) node[midway, below] {$k_{10}$};
\draw [edge]  (9) to node[above] {$k_{11}$} (10);

\draw [edge] ([yshift=2pt]11.east) -- ([yshift=2pt]12.west) node[midway, above] {$k_{12}$};
\draw [edge] ([yshift=-2pt]12.west) -- ([yshift=-2pt]11.east) node[midway, below] {$k_{13}$};
\draw [edge]  (12) to node[above] {$k_{14}$} (13);

\end{tikzpicture}
\end{center}

This can be decomposed into two independent subnetworks. The first subnetwork is given by
\begin{center}
\begin{tikzpicture}[baseline=(current  bounding  box.center)]
\tikzset{vertex/.style = {draw=none,fill=none}}
\tikzset{edge/.style = {->,> = latex', line width=0.10mm}} 
\node[vertex] (1) at  (-4,0) {$XD$};
\node[vertex] (2) at  (-2,0) {$X$};

\draw [edge] ([yshift=2pt]1.east) -- ([yshift=2pt]2.west) node[midway, above] {$k_1$};
\draw [edge] ([yshift=-2pt]2.west) -- ([yshift=-2pt]1.east) node[midway, below] {$k_2$};
\end{tikzpicture}
\end{center}
and is already WR and DZ. The second subnetwork, however, requires translation. It is shown below together with the translation complexes to achieve WR and DZ. 
\begin{center}
\begin{tikzpicture}[baseline=(current  bounding  box.center)]
\tikzset{vertex/.style = {draw=none,fill=none}}
\tikzset{edge/.style = {->,> = latex', line width=0.10mm}} 
\node[vertex] (2) at  (-2,0) {$X$};
\node[vertex] (3) at  (0,0) {$XT$};
\node[vertex] (4) at  (2,0) {$X_p$};

\node[vertex] (14) at  (5,0) {$(+XD+XT+Y_p)$};

\node[vertex] (5) at  (-2,-1.5) {$X_p + Y$};
\node[vertex] (6) at  (0,-1.5) {$X_p Y$};
\node[vertex] (7) at  (2,-1.5) {$X + Y_p$};

\node[vertex] (15) at  (5,-1.5) {$(+XD+XT)$};

\node[vertex] (8) at  (-2,-3) {$XD + Y_p$};
\node[vertex] (9) at  (0,-3) {$XDY_p$};
\node[vertex] (10) at  (2,-3) {$XD + Y$};

\node[vertex] (16) at  (5,-3) {$(+XT + X_p)$};

\node[vertex] (11) at  (-2,-4.5) {$XT + Y_p$};
\node[vertex] (12) at  (-0,-4.5) {$XTY_p$};
\node[vertex] (13) at  (2,-4.5) {$XT + Y$};

\node[vertex] (17) at  (5,-4.5) {$(+XD + X_p)$};

\draw [edge] ([yshift=2pt]2.east) -- ([yshift=2pt]3.west) node[midway, above] {$k_3$};
\draw [edge] ([yshift=-2pt]3.west) -- ([yshift=-2pt]2.east) node[midway, below] {$k_4$};
\draw [edge]  (3) to node[above] {$k_5$} (4);

\draw [edge] ([yshift=2pt]5.east) -- ([yshift=2pt]6.west) node[midway, above] {$k_6$};
\draw [edge] ([yshift=-2pt]6.west) -- ([yshift=-2pt]5.east) node[midway, below] {$k_7$};
\draw [edge]  (6) to node[above] {$k_8$} (7);

\draw [edge] ([yshift=2pt]8.east) -- ([yshift=2pt]9.west) node[midway, above] {$k_9$};
\draw [edge] ([yshift=-2pt]9.west) -- ([yshift=-2pt]8.east) node[midway, below] {$k_{10}$};
\draw [edge]  (9) to node[above] {$k_{11}$} (10);

\draw [edge] ([yshift=2pt]11.east) -- ([yshift=2pt]12.west) node[midway, above] {$k_{12}$};
\draw [edge] ([yshift=-2pt]12.west) -- ([yshift=-2pt]11.east) node[midway, below] {$k_{13}$};
\draw [edge]  (12) to node[above] {$k_{14}$} (13);

\end{tikzpicture}
\end{center}    

By applying CRITERIA, we obtain the following equilibrium parametrization
\begin{align*}
    X &= \dfrac{X_pY k_8 (k_4+k_5)}{k_3 k_5} \\
    XD &= \dfrac{X_pY k_2 k_8 (k_4+k_5)}{k_1 k_3 k_5}  \\
    XDY_p &= \dfrac{X_pY k_8 k_9 (k_{13} + k_{14})}{k_9 k_{11} k_{13} + k_9 k_{11} k_{14} + \dfrac{k_1 k_3 k_{10} k_{12} k_{14}}{k_2 (k_4 + k_5)} + \dfrac{k_1 k_3 k_{11} k_{12} k_{14}}{k_2 (k_4 + k_5)}}  \\
    XT &= \dfrac{X_pY k_8}{k_5} \\
    XTY_p &= \dfrac{X_pY k_1 k_3 k_8 k_{12} (k_{10} + k_{11})}{k_2 (k_4 + k_5) \left( k_9 k_{11} k_{13} + k_9 k_{11} k_{14} + \dfrac{k_1 k_3 k_{10} k_{12} k_{14}}{k_2 (k_4 + k_5)} + \dfrac{k_1 k_3 k_{11} k_{12} k_{14}}{k_2 (k_4 + k_5)} \right)} \\
    X_p &= \dfrac{X_pY (k_7 + k_8)}{Y k_6} \\
    Y_p &= \dfrac{k_1 k_3 k_5 (k_{10} + k_{11}) (k_{13} + k_{14})}{k_2 (k_4 + k_5) \left( k_9 k_{11} k_{13} + k_9 k_{11} k_{14} + \dfrac{k_1 k_3 k_{10} k_{12} k_{14}}{k_2 (k_4 + k_5)} + \dfrac{k_1 k_3 k_{11} k_{12} k_{14}}{k_2 (k_4 + k_5)} \right)}
\end{align*}
with free parameters $X_pY$ and $Y$, and set of conservation laws
\begin{align*}
    \dfrac{d X}{dt} + \dfrac{d XD}{dt} + \dfrac{d XDY_p}{dt} + \dfrac{d XT}{dt} + \dfrac{d XTY_p}{dt} + \dfrac{d X_p}{dt} + \dfrac{d X_pY}{dt} &= 0 \\
    -\dfrac{d X}{dt} - \dfrac{d XD}{dt} - \dfrac{d XT}{dt} - \dfrac{d X_p}{dt} + \dfrac{d Y}{dt} + \dfrac{d Y_p}{dt} &= 0
\end{align*}

Observe that the equilibrium concentration of $Y_p$ depends only on rate constants confirming the result in \cite{JMP2019:parametrization}. Now, integrating the conservation equations with respect to $t$, we get
\begin{align}
    X + XD + XDY_p + XT + XTY_p + X_p + X_pY &= T_1 \label{eq:T1} \\
    -X-XD-XT-X_p + Y + Y_p &= T_2 \label{eq:T2}
\end{align}
where $T_1$ and $T_2$ are nonnegative constants dependent on the initial conditions of the species involved in their respective equation. To facilitate our computation, we express the rate constants in the equilibrium parametrization in terms of alphas so that $X = X_pY \alpha_1$, $XD = X_pY \alpha_2$, $XDY_p = X_pY \alpha_3$, $XT = X_pY \alpha_4$, $XTY_p = X_pY \alpha_5$, $X_p = \dfrac{X_pY \alpha_6}{Y}$, and $Y_p = \alpha_7$, where $\alpha_i > 0$ for $i=1,2,3,4,5,6,7$. Plugging these expressions to Eq \ref{eq:T1} and Eq \ref{eq:T2}, we get
\begin{align}
    X_pY\alpha_1 + X_pY\alpha_2 + X_pY\alpha_3 + X_pY\alpha_4 + X_pY\alpha_5 + \dfrac{X_pY\alpha_6}{Y} + X_pY &= T_1 \label{eq:T1_2} \\
    -X_pY\alpha_1 - X_pY\alpha_2 - X_pY\alpha_3 - \dfrac{X_pY\alpha_6}{Y} + Y + \alpha_7 &= T_2 \label{eq:T2_2}
\end{align}

From Eq \ref{eq:T1_2}, we have $Y = \dfrac{X_pY \alpha_6}{T_1 - X_pY \gamma_1}$, where $\gamma_1 = 1 + \alpha_1 + \alpha_2 + \alpha_3 + \alpha_4 + \alpha_5$. Similarly, we can write Eq \ref{eq:T2_2} as
\begin{align*}
    -Y \cdot X_pY \gamma_2 - X_pY \alpha_6 + Y^2 + Y\alpha_7 = Y \cdot T_2
\end{align*}
where $\gamma_2 = \alpha_1 + \alpha_2 + \alpha_3$ so that
{\scriptsize
\begin{align*}
    -\left( \dfrac{X_pY \alpha_6}{T_1 - X_pY \gamma_1} \right) X_pY \gamma_2 - X_pY \alpha_6 + \left( \dfrac{X_pY \alpha_6}{T_1 - X_pY \gamma_1} \right)^2 + \left( \dfrac{X_pY \alpha_6}{T_1 - X_pY \gamma_1} \right)\alpha_7 &= \left( \dfrac{X_pY \alpha_6}{T_1 - X_pY \gamma_1} \right) T_2 \\
    \implies \left(\gamma_1 \alpha_6 \gamma_2 - \gamma_1^2 \alpha_6\right) X_pY^3 + \left(2T_1 \gamma_1 \alpha_6 - T_1\alpha_6 \gamma_2 + \alpha_6^2 - \gamma_1 \alpha_6\alpha_7 + \gamma_1 \alpha_6 T_2\right)X_pY^2 &+ \\ (T_1\alpha_6\alpha_7 - T_1^2\alpha_6 - T_1\alpha_6 T_2)X_pY &= 0
\end{align*}
}

\begin{figure*}[h!]
\centering
\includegraphics[width=13.5cm,height=15.5cm,keepaspectratio]{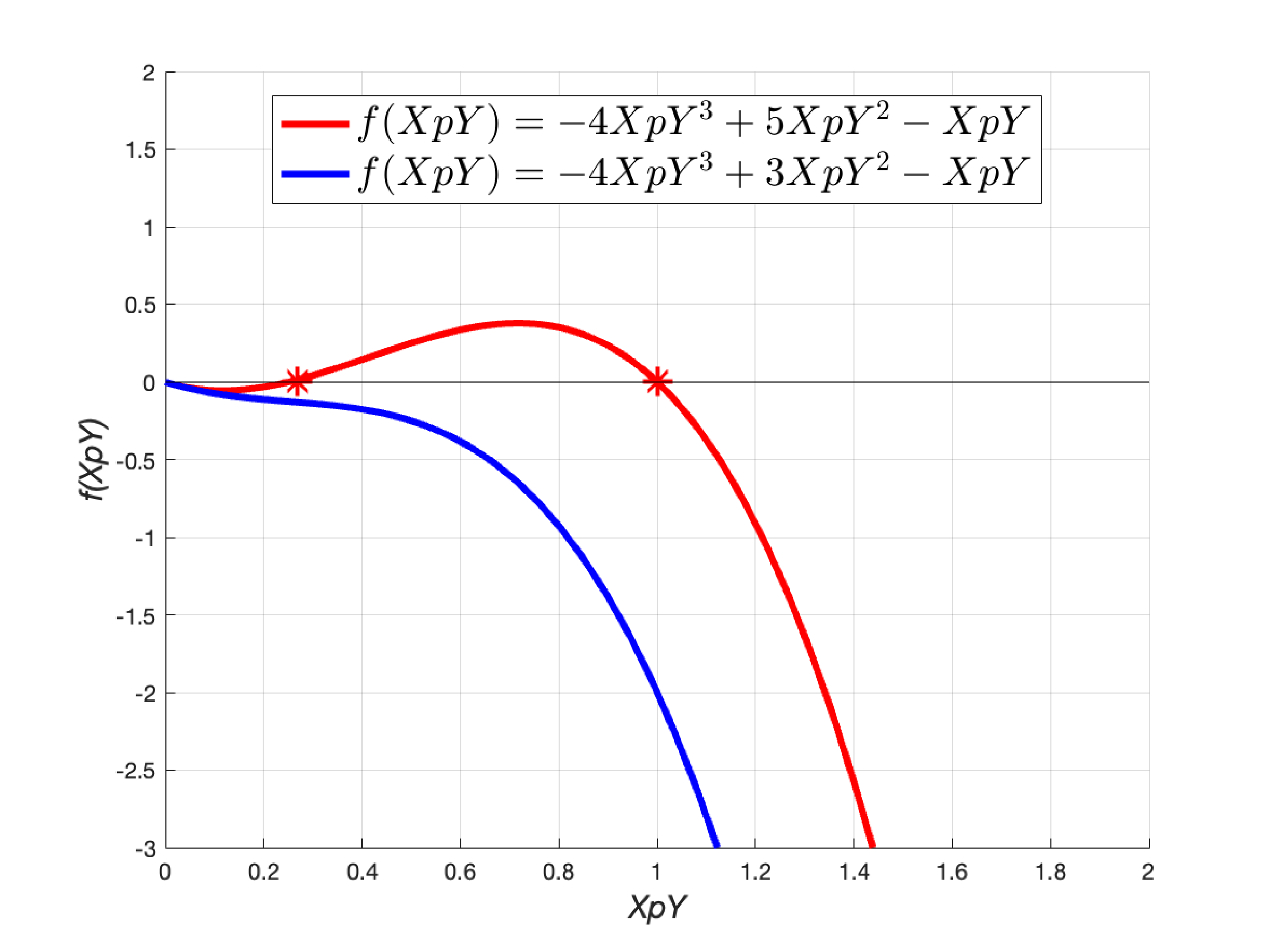}
\caption{\scriptsize
{\bf{Graphs of two functions of the form $\mathbf{f(X_pY) = -a_3X_pY^3 + a_2X_pY^2 - a_1X_pY}$ with positive constants $\mathbf{a_i}$.}} The $x$-axis represents the value of $X_pY$ at steady state. The first function, shown in red, has two positive zeroes while the second function, shown in blue, has no positive zero.
}\label{fig:multistationarity}
\end{figure*}

The left hand side of the equation is a cubic polynomial in $X_pY$ with two sign changes indicating that the equation can only have two or zero positive real solutions by Descartes' Rule of Signs. Refer to the Supplementary Information for the proof. Equivalently, $X_pY$ can admit exactly two positive equilibria or none at all, implying that the system has the capacity for multistationarity. This behavior is supported by Fig \ref{fig:multistationarity} which shows representative cases for specific choices of rate constants and initial conditions. In particular, the first function exhibits two positive roots, whereas the second admits none.

In the context of the EnvZ-OmpR system, the existence of two positive equilibria or none at all provides a mechanistic explanation for how the EnvZ–OmpR pathway can implement threshold-based osmotic sensing. One equilibrium may represent basal regulation, while the second may correspond to a fully activated stress-response state. This provides a rigorous mathematical foundation for its role as a master regulator of osmotic stress response. This result demonstrates the predictive power of CRITERIA over numerical simulations. While the latter can show specific cases of multistationarity, the analytical result provides a guarantee that no more than two positive steady states can exist regardless of the chosen rate constants. CRITERIA reveals how the network's structure, rather than just specific parameter values, dictates these functions.

For our next application, we consider the full CRISPRi toggle switch model which consists of $17$ species and $42$ reactions \cite{SantosMoreno2020}. This is Model Z of Table \ref{tab:networkproperties}. In this system, a catalytically inactive protein (dCas9) forms complexes with two single-guide RNAs (sg1 and sg2). These complexes subsequently bind both to their target genes ($G_1$ and $G_2$, respectively) and also to the opposing genes \cite{SantosMoreno2020,HernandezPLOS2023}. Binding to the intended target genes is referred to as specific binding, whereas binding to the opposing genes is called unspecific binding. The analysis in \cite{SantosMoreno2020} shows that the interplay between specific and unspecific binding of dCas9/sgRNA complexes to the genes explains bistability in the CRISPRi toggle switch model.

\begin{figure*}[h!]
\centering
\includegraphics[width=13.5cm,height=35.5cm,keepaspectratio]{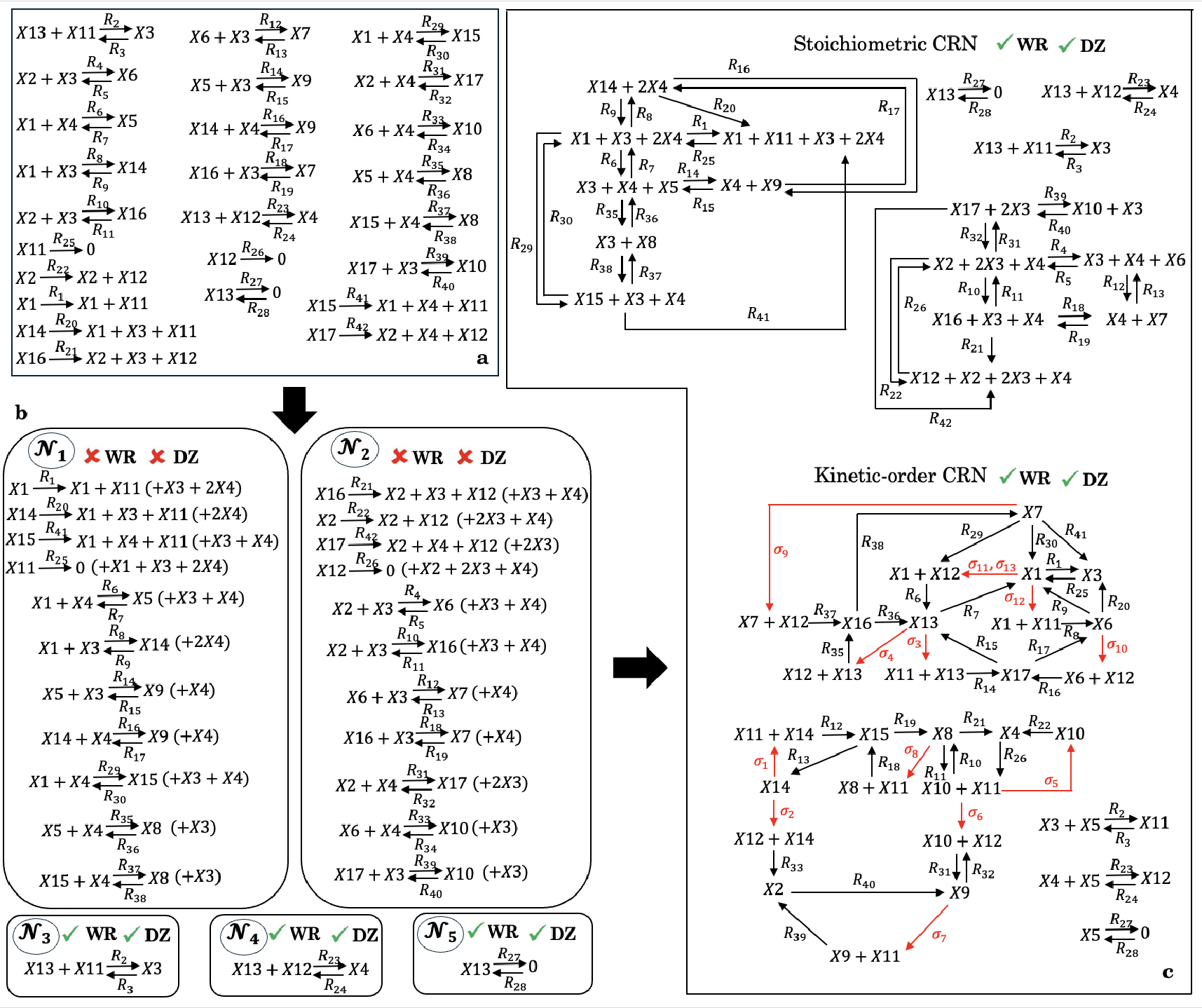}
\caption{\scriptsize
{\bf{Network decomposition, translation, and merging of the full CRIPSRi model for equilibrium parametrization.}} 
{\bf{a}} The reaction network consists of $17$ species and $42$ reactions. The species of interest in our analysis is $X13$, which represents the catalytically inactive dCas9 protein. The full correspondence between model species and network variables is provided in the Supplementary Information.
{\bf{b}} The network has five independent subnetworks, two of which are neither WR nor DZ. The other three subnetworks consist only of two reactions which are reversible. The complexes to the right of each reaction are the translation complexes after applying the elementary flux mode-based translation method described above. 
{\bf{c}} Adding the translation complexes to each reaction yields the following stoichiometric CRN and kinetic-order CRN, which together constitute the generalized reaction network for the full CRISPRi model. The stoichiometric CRN is the resulting network after adding the translation complexes to each of the reactions in the original model. On the other hand, the kinetic-order CRN, which determines the kinetics of the system, consists of the original source complex of the reaction network before translation, thereby ensuring that the dynamics is preserved. The reactions shown in red in the kinetic-order CRN are referred to as phantom edges. They connect identical complexes in the stoichiometric CRN. In the example, both the stoichiometric and kinetic-order CRNs are WR and DZ. Hence, the translated network is also WR and DZ. This now enables the parametrization of the equilibria of the model.
}\label{fig:crispri}
\end{figure*}

Here, we show that, in addition to the system's capacity for multistationarity, the model also exhibits ACR in the dCas9 species. To establish this, we parametrize the equilibria of the system by applying CRITERIA. Fig \ref{fig:crispri} illustrates the decomposition, translation, and subsequent merging of the reaction network prior to deriving the equilibrium expressions. The species of interest in this analysis is $X13$, which represents the dCas9 protein. As shown in Fig \ref{fig:crispri}b, the network has five independent subnetworks, two of which are neither WR nor DZ. Using the elementary flux mode-based network translation outlined above, we apply the translation complexes indicated to the right in each of the reactions in Fig \ref{fig:crispri}b produding a generalized reaction network that is both WR and DZ. 

A generalized reaction network consists of two components: the \textit{stoichiometric CRN}, which contains the complexes obtained after adding the translation complexes to the reactions, and the \textit{kinetic-order CRN}, which retains the original source complexes of the reaction network, thereby preserving the dynamics of the system. The reactions in red in the kinetic-order CRN are called \textit{phantom edges}, which connect identical complexes in the stoichiometric CRN. If both stoichiometric CRN and kinetic-order CRN are WR and DZ, then the translated network is WR and DZ. See Methods section for more information.

Using CRITERIA, the derived equilibrium expression for dCas9 is given by: $$\text{X13} = \frac{k_{28}}{k_{27}}$$ where $k_{28}$ and $k_{27}$ represent the rate constants for the production ($R_{28}$) and degradation ($R_{27}$) of dCas9, respectively. This tells us that the equilibrium concentration of the dCas9 species is dependent only on these specific rate constants. So, the equilibrium concentration of dCas9 remains robustly constant for any initial conditions under fixed parameters. This remarkable property is known as \textit{absolute concentration robustness} (ACR) \cite{ShinarFeinberg2010}. 
Such finding has implications for the design of robust synthetic systems:
\begin{itemize}
\item \textbf{Buffering Against Network Fluctuations}: ACR suggests an intrinsic buffering mechanism wherein the equlibrium concentration of dCas9 remains invariant despite fluctuations in other network components, such as target gene concentrations or guide RNA levels.
\item \textbf{Predictable Tuning}: In synthetic biology, achieving predictable and robust behavior is a primary concern. Our discovery shows that the long-term concentration of dCas9 is robustly constant over the concentrations of all other species.
\item \textbf{Simplified Design Parameters}: Fine-tuning the long-term concentration of dCas9 in the circuit can only be achieved by altering its production rate $k_{28}$ or degradation rate $k_{27}$. This decoupling from other parameters allows for modular circuit design where dCas9 levels can be set independently of the broader network complexity.
\end{itemize}
By maintaining a stable pool of the dCas9 engine, the CRISPRi circuit ensures that the machinery of regulation remains consistent across varying cellular environments, effectively anchoring the stability of the toggle switch.

\subsection*{Computational package, CRITERIA}

We developed CRITERIA (https://github.com/evvillejo/CRITERIA), a user-friendly, open-source, and publicly available computational package that automates the computation of equilibria of biochemical reaction networks. The framework enhances and extends COMPILES by incorporating several methodological improvements that overcome previously identified limitations (see Fig \ref{fig:CRITERIAmethod} and its accompanied discussion).

Taking the set of reactions of a CRN as input, CRITERIA outputs the following: an equilibrium parametrization, kinetic deficiency of the GCRN along with the additional equations required by Theorem \ref{thm:deficiency-based}, and the conservation laws of the system. As with COMPILES, CRITERIA does its best to output the most simplified analytic solution in terms of rate constants and free parameters.


\section*{Discussion}

\begin{figure*}
\centering
\includegraphics[width=13.5cm,height=15.5cm,keepaspectratio]{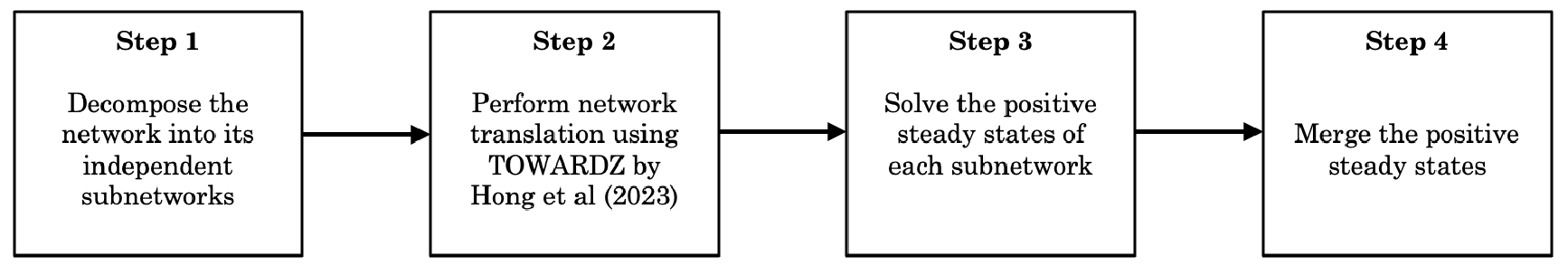}
\caption{\scriptsize
{\bf{General framework of COMPILES for deriving analytic equilibria.}}
The method starts by decomposing the network into its independent subnetworks, thereby simplifying the analysis by focusing on smaller individual pieces rather than looking the network as a whole. Next, any subnetwork that is not WR and DZ is translated using the TOWARDZ approach introduced in \cite{HongSIAM}. Once all subnetworks satisfy WR and DZ, the procedure proceeds by determining the positive equilibria of each subnetwork through the parametrization method in \cite{JMP2019:parametrization}. Finally, the equilibria of the original network are obtained by merging the equilibria of the subnetworks. 
}\label{fig:COMPILESmethod}
\end{figure*}

In this paper, we have developed an enhanced framework together with a corresponding computational package for computing the positive equilibria of biochemical reaction networks. This approach substantially extends and improves the method introduced in \cite{HernandezPLOS2023}, called COMPILES, for several key reasons (see Fig \ref{fig:COMPILESmethod} for the schematic diagram of its framework):
\begin{itemize}
    \item Its network translation step constitutes a major computational bottleneck in the procedure. In particular, TOWARDZ, the package used to perform network translation in COMPILES, takes too long to output a WR and DZ translation.
    \item It has a ``merging'' issue that leads to interdependencies in the resulting equilibrium parametrization.    
    \item COMPILES is limited only to reaction networks with zero kinetic deficiency.
\end{itemize}

\begin{figure}[H]
\centering
\includegraphics[width=\linewidth,keepaspectratio]{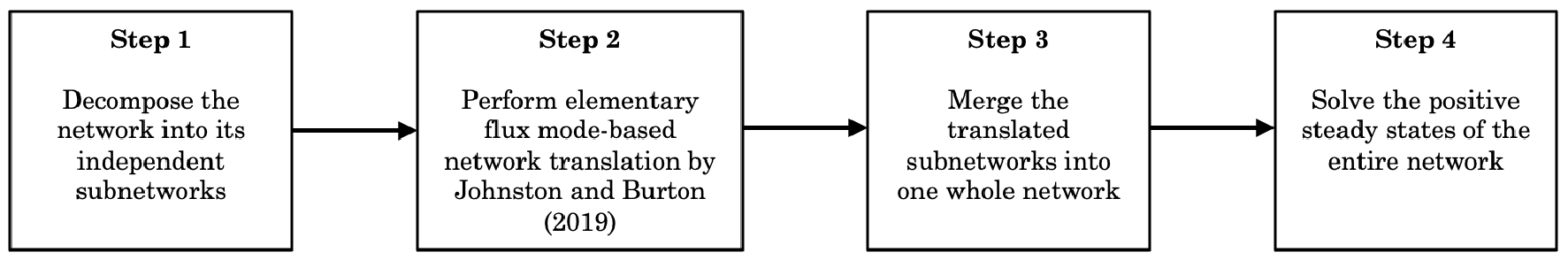}
\caption{\scriptsize
{\bf{An enhanced framework CRITERIA for deriving analytic equilibria.}}
Similar to COMPILES, the framework for CRITERIA begins with independent network decomposition which is then followed by network translation for each independent subnetwork that is not WR and DZ using the elementary flux mode-based approach \cite{JohnstonBurton2019}. In contrast to COMPILES, however, CRITERIA first performs merging of the translated subnetworks into a single network by taking the union of the reactions in the subnetworks before proceeding to compute the positive equilibria of the entire network, which constitutes the last step of the enhanced procedure.
}
\label{fig:CRITERIAmethod}
\end{figure}

Our approach addresses the issues and limitations enumerated above by adopting the elementary flux mode-based translation approach in \cite{JohnstonBurton2019}, reordering the steps in the overall procedure of COMPILES (i.e., we first merge the translated subnetworks into one whole network before solving the equilibria), and incorporating additional conditions required to parametrize the equilibria of networks with positive kinetic deficiency (Theorem 15 of \cite{JMP2019:parametrization}). These modifications form the foundation of the computational package we developed, termed CRITERIA, which implements the improved process. The general framework of CRITERIA is shown in Fig \ref{fig:CRITERIAmethod}.

We were able to show that CRITERIA indeed outperforms COMPILES, particularly for large reaction networks, in terms of computational runtime. (see Fig \ref{tab:towardzVSrewardz}). In fact, COMPILES was unable to parametrize the equilibria of 11 out of the 26 models considered in this study. Moreover, 10 of the 26 models have positive kinetic deficiency, highlighting a class of networks that COMPILES cannot accommodate but are handled by CRITERIA. 

It has been observed that COMPILES may handle smaller networks slightly more efficiently than CRITERIA, as the latter requires solving elementary flux modes and a binary linear program, both of which are not trivial tasks. This advantage arises when independent decomposition directly yields WR and DZ subnetworks, thereby eliminating the need for translation. However, such favorable decompositions cannot be expected in every biochemical system, as they depend strongly on the network's structure and interconnections. Indeed, the runtime differences observed for the 2-Protein Gene Transcription Model, Insulin Signaling Pathway, and CRISPRi Model 1 contrast sharply with those for the Hybrid Histidine–Kinase Model, MAPK Model, and Zigzag Model, where CRITERIA significantly outperformed COMPILES.

The improvements in the computational runtime are largely attributable to the integration of the elementary flux mode-based network translation, which removes a key bottleneck in COMPILES. Compared with the earlier approach, this translation strategy is more structured and systematic, leading to enhanced numerical stability and improved scalability. Furthermore, the observations above suggest that the improved ordering of steps in CRITERIA alleviates the computational burden arising from handling subnetworks separately, thereby improving the parametrization stage of the procedure. As a result, the interdependencies that previously appeared in the derived equilibrium expressions are effectively resolved. Taken together, the results indicate that CRITERIA is not merely a computational enhancement but also a conceptual extension of the original framework. It broadens the range of networks for which analytic equilibria can be derived while maintaining practical efficiency, making it a promising tool for the systematic analysis of increasingly large and structurally complex biochemical reaction networks.

Finally, we demonstrate the utility of having a closed form equilibrium expression for analyzing key dynamical properties of biochemical systems. Using CRITERIA, we were able to easily obtain the analytic equilibria of the EnvZ-OmpR model which COMPILES was not able to do within the prescribed maximum runtime. Using the system's conservation laws, we showed that the model has the capacity to admit multiple positive steady states. In addition to exhibiting absolute concentration robustness, as previously established in \cite{ShinarFeinberg2010,JMP2019:parametrization}, our analysis provides further insight that the system admits exactly two positive equilibria or none at all. As a second application, we easily obtained an equilibrium parametrization of the full CRISPRi toggle switch model, which COMPILES also failed to compute within the prescribed maximum runtime. Our analysis reveals that the dCas9 species exhibits ACR, indicating that its concentration is invariant across all positive equilibria and independent of initial conditions. This robustness property suggests a potential buffering mechanism in CRISPRi-based regulatory circuits, by which the availability of dCas9 is maintained despite fluctuations elsewhere in the network. From a design perspective, this feature may be advantageous in synthetic biology settings, where predictable regulator levels are essential for reliable circuit performance.

These applications demonstrate how explicit equilibrium parametrizations enable a more general analysis of biochemical systems beyond simulation-based methods. In particular, closed-form equilibrium expressions make analysis of models possible without committing to specific parameter values, thereby separating intrinsic dynamical features from numerical artifacts. By enabling analytic equilibrium parametrizations for broader classes of biochemical reaction networks, CRITERIA opens the door to parameter-free analyses of complex models, allowing one to characterize multistationarity, robustness, and other qualitative dynamical behavior directly from network structure. This capability not only reduces dependence on uncertain or experimentally inaccessible parameter estimates but also enables systematic structural comparisons between models and reveals general mechanistic principles underlying their dynamical behavior. Looking ahead, this study presents several avenues for future work.
\begin{enumerate}
    \item The translation method only works for elementary flux modes that are unitary (see Theorem \ref{thm: structuraltranslation} in the Methods section). Future work will focus on adapting the translation procedure to account for models whose elementary flux modes are not unitary.
    \item A promising direction for future work is to assess multistationarity using Sturm's Theorem, which is a result in algebraic geometry \cite{SiegalGaskins2015}. This approach follows a systematic sequence of steps that ultimately yields a polynomial in the species concentrations, with coefficients depending on the system's parameters. By analyzing the sign structure of this polynomial, one can infer the existence or absence of multiple equilibria as a function of both parameters and variables. A key prerequisite for this method is an equilibrium parametrization, which our current framework already provides. As a next step, we plan to automate this method by integrating it with the equilibrium parametrization generated by our approach.
    \item Another important direction for further study is the automation of boundary equilibrium analysis, where the equilibrium value of one or more species is zero. This approach would extend the framework to include biologically relevant cases such as extinction states or other switch-like behaviors, offering a more comprehensive understanding of the system's dynamics.
\end{enumerate}

More broadly, building on the current framework could support large-scale, structure-driven analyses of complex biochemical networks, unlocking new opportunities to translate theoretical insights from reaction network theory into practical tools for systems and synthetic biology.


\section*{Methods}

\subsection*{Chemical reaction networks}

A \textit{chemical reaction network} (CRN) is defined by three sets \cite{FeinbergBook2019}. The first set is the set of species $\mathcal{S}= \{X_1,\ldots,X_m \}$, which contains the fundamental units of the network called \textit{species}. The second set is the set of complexes $\mathcal{C}=\{y_1,\ldots,y_n\}$, which contains entities called \textit{complexes} that are linear combinations of species. That is,
\begin{align*}
    y_i = \sum\limits_{j=1}^m y_{ij} X_j, \hspace{.3cm} i=1,\ldots,n.
\end{align*}

The coefficients $y_{ij} \in \mathbb{Z}_{\geq0}$ are called stoichiometric coefficients. We use $y_i$ as both the complex itself and the corresponding complex vector $y_i = (y_{i1},\ldots,y_{im}) \in \mathbb{Z}_{\geq 0}^m$. And lastly, the third set is the set of reactions $\mathcal{R} = \{r_1,\ldots,r_r\} \subset \mathcal{C} \times \mathcal{C}$.

The structure of a CRN can be seen as a directed graph where the edges are the reactions, and the vertices are the complexes. So, we can represent reactions as ordered pairs of complexes (i.e., $r_k = \left( y_i,y_j \right)$) or as directed edges (i.e., $r_k = y_i \to y_j$). Here, we call $y_i$ as the \textit{reactant} or \textit{source complex} and $y_j$ as the \textit{product complex}. A weight $k_i \in \mathbb{R}_{\geq0}$ is associated to each reaction called the \textit{reaction rate coefficient}. For each reaction $r_k = y_i \to y_j$, we can associate a \textit{reaction vector} $y_j - y_i \in \mathbb{Z}^m$ which serves as the columns of the \textit{stoichiometric matrix} $N$ of the CRN.

\begin{example}
\label{main:example:CRN}
Consider the following toy CRN with five reactions:
\begin{center}
\begin{tikzpicture}[baseline=(current bounding box.center)]
\tikzset{vertex/.style = {draw=none,fill=none}}
\tikzset{edge/.style = {->, >=latex', line width=0.10mm}} 

\node[vertex] (1) at  (-1,0) {$R_1:A+B$};
\node[vertex] (2) at  (1,0) {$C+B$};
\node[vertex] (3) at  (3,0) {$R_2,R_3:A$};
\node[vertex] (4) at  (5,0) {$B$};
\node[vertex] (5) at  (6.5,0) {$R_4:C$};
\node[vertex] (6) at  (8,0) {$0$};
\node[vertex] (7) at  (9,0) {$R_5:0$};
\node[vertex] (8) at  (10.5,0) {$A$};

\draw[edge]  (1) -- node[above] {$k_1$} (2);

\draw [edge] ([yshift=2pt]3.east) -- ([yshift=2pt]4.west) node[midway, above] {$k_2$};
\draw [edge] ([yshift=-2pt]4.west) -- ([yshift=-2pt]3.east) node[midway, below] {$k_3$};

\draw[edge]  (5) -- node[above] {$k_4$} (6);
\draw[edge]  (7) -- node[above] {$k_5$} (8);
\end{tikzpicture}
\end{center}

In this CRN, $R_1$ represents the conversion of $A$ to $C$ in the presence of $B$, with $B$ effectively acting as a catalyst for the reaction. The reversible reactions $R_2$ and $R_3$ can be interpreted as conversion of $A$ to $B$, and vice versa. Reaction $R_4$ models the degradation or removal of $C$ from the system, while $R_5$ represents the inflow or external production of $A$. The network consists of three species, six complexes, and five reactions. In particular,
\begin{align*}
    \mathcal{S} &= \{ A,B,C \}, \\
    \mathcal{C} &= \{ A+B,C+B,C,A,B,0 \}, \\
    \mathcal{R} &= \{ A+B \to C+B, A \to B, B \to A,C\to 0, 0 \to C \}.
\end{align*}

The zero complex is not considered as a species; rather, it is a formal placeholder that represents the external environment of the system, where it acts as a "sink" or "source" of species and complexes. In $R_1 : A+B \to C+B$, $A+B$ acts as the source complex and $C+B$ acts as the product complex. Furthermore, the following are the reaction vectors of the CRN:
\begin{align*}
    R_1 &: (C+B) - (A+B) = C-A=[-1 \hspace{.2cm} 0 \hspace{0.2cm} 1]^T \hspace{1cm} R_4 : (0) - (C) = -C = [0 \hspace{.2cm} 0 \hspace{0.2cm} -1]^T \\
    R_2 &: (B) - (A) = B-A = [-1 \hspace{.2cm} 1 \hspace{0.2cm} 0]^T \hspace{2.4cm} R_5 : (A) - (0) = A = [1 \hspace{.2cm} 0 \hspace{0.2cm} 0]^T  \\
    R_3 &: (A) - (B) = A-B = [1 \hspace{.2cm} -1 \hspace{0.2cm} 0]^T
\end{align*}

These reaction vectors comprise the columns of the stoichiometric matrix, which can be written as 
\begin{equation*}
    N=\begin{blockarray}{cccccc}
    R_1 & R_2 & R_3 & R_4 & R_5 \\
    \begin{block}{[rrrrr]c}
    -1 & -1 & 1 & 0 & 1 & A \\ 
    0 & 1 & -1 & 0 & 0 & B \\
    1 & 0 & 0 & -1 & 0 & C \\
\end{block}
\end{blockarray}
\end{equation*}
whose rank is given by $s=3$.
\end{example}

Since a CRN can be seen as a directed graph, we can draw parallels between concepts in graph theory and CRN theory. In graph theory, we say that two vertices are connected if there is a directed path between them. Furthermore, two vertices are strongly connected if there is a directed path to and from each other. Given these definitions, we call the connected components of the network as the \textit{linkage classes} of the CRN and the strongly connected components as the \textit{strong linkage classes} of the CRN. In other words, a linkage class is a set of complexes that are connected to one another but not to any other complex not in the set, whereas a strong linkage class consists of complexes that are strongly connected to each other. If the number of linkage classes is equal to the number of strong linkage classes, then the CRN is said to be \textit{weakly reversible}. Alternatively, a CRN is said to be \textit{weakly reversible} if every reaction is part of a directed cycle.

\begin{example}
In Example \ref{main:example:CRN}, there are two linkage classes, namely $\{A+B,C+B\}$ and $\{ A,B,C,0 \}$. On the other hand, there are five strong linkage classes which are the following: $\{ A,B \}$, $\{ C+B \}$, $\{ A+B \}$, $\{ C \}$, $\{ 0 \}$. Note that a complex, if not strongly connected to other complex/es, is a trivial strong linkage class. Since the number of linkage classes does not coincide with the number of strong linkage classes, the CRN is not weakly reversible. Alternatively, this can be seen from the fact that reaction $R_1$ is not contained in any directed cycle, and hence the network fails to be weakly reversible.
\end{example}

One of the central concepts in CRN theory is called the \textit{deficiency} of the CRN which is denoted by $\delta$. It is a nonnegative integer with formula $\delta = n - l - s$, where $n$ is the number of complexes, $l$ is the number of linkage classes of the CRN, and $s$ is the rank of the stoichiometric matrix. It can be interpreted as a measure of linear dependence among the reactions of the reaction network. The higher the deficiency, the lower the extent of linear independence \cite{ShinarFeinberg2011}.

\begin{example}
    In Example \ref{main:example:CRN}, we have $\delta=6-2-3=1$. So, it is a deficiency one network.
\end{example}

A CRN is endowed by kinetics to describe the evolution of the concentration of the species through time. One of the widely-used kinetics is \textit{mass-action} which states that the rate at which the reaction proceeds is proportional to the product of the concentration of the species in its reactant complex.

\begin{example}
    In our running example, suppose $a$, $b$, and $c$ denote the concentrations of species $A, B,$ and $C$, respectively. Then, the rate function of each reaction under the assumption of mass action is given by
\begin{equation*}
    R_1: k_1 ab \hspace{.5cm} R_2: k_2 a \hspace{.5cm} R_3: k_3 b \hspace{.5cm} R_4: k_4 c \hspace{.5cm} R_5: k_5.
\end{equation*}

Together with $N$, the system of ordinary differential equations for the model is given by
\begin{equation*}
    \begin{bmatrix}
        \dfrac{da}{dt} \\[1em]
        \dfrac{db}{dt} \\[1em]
        \dfrac{dc}{dt} \\
    \end{bmatrix} =
    k_1 ab
    \begin{bmatrix*}[r]
        -1 \\
        0 \\
        1 
    \end{bmatrix*} + 
    k_2a
    \begin{bmatrix*}[r]
        -1 \\
        1 \\
        0 
    \end{bmatrix*} + 
    k_3 b
    \begin{bmatrix*}[r]
        1 \\
        -1 \\
        0
    \end{bmatrix*} +
    k_4 c
    \begin{bmatrix*}[r]
        0 \\
        0 \\
        -1
    \end{bmatrix*} + 
    k_5
    \begin{bmatrix*}[r]
        1 \\
        0 \\
        0
    \end{bmatrix*},
\end{equation*}
where $\dfrac{da}{dt}$, $\dfrac{db}{dt}$, and $\dfrac{dc}{dt}$ are the time derivatives of the concentration functions of the species $A$, $B$, and $C$, respectively. The system of ODEs can be written as
\begin{align*}
    \dfrac{da}{dt} &= -k_1 ab - k_2 a + k_3 b + k_5 \\
    \dfrac{db}{dt} &= k_2a - k_3b \hspace{.2cm} \cdot \\
    \dfrac{dc}{dt} &= k_1 ab - k_4 c
\end{align*}
\end{example}

The long-term behavior of a system is described by the \textit{equilibrium} or \textit{steady states}, which can be solved by setting the time derivatives to zero and then solving for the concentration of the species.

\begin{example}
    The positive equilibrium of the network in Example \ref{main:example:CRN} is given by
    \begin{align*}
        \left( a,b,c \right) = \left( \sqrt{\dfrac{k_3 k_5}{k_1 k_4}}, \sqrt{\dfrac{k_3 k_4}{k_1 k_5}}, \dfrac{k_3}{k_2} \right),
    \end{align*}
    where $k_1,k_2,k_3,k_4,k_5 > 0$.
\end{example}

\subsection*{Network decomposition}

A \textit{decomposition} of a CRN is induced by a partition of its reaction set $\mathcal{R}$ \cite{HernandezPLOS2023}. This partition then gives rise to a decomposition of the reaction network into subnetworks. One special type of network decomposition is called \textit{independent decomposition}, where the rank of the stoichiometric matrix of the whole network is equal to the sum of the ranks of the stoichiometric matrices of its subnetworks \cite{FeinbergBook2019}.

\begin{example}
    In Example \ref{main:example:CRN}, the following partition of the reaction set $\mathcal{R}$ into
\begin{align*}
    R_1 : \{ R_1, R_4, R_5 \} \hspace{0.5cm} \text{and} \hspace{0.5cm} R_2 : \{ R_2, R_3 \}
\end{align*}
induces a decomposition of the CRN into subnetworks $\mathcal{N}_1$ and $\mathcal{N}_2$ with stoichiometric matrices
\begin{equation*}
    \begin{blockarray}{cccc}
    & R_1 & R_4 & R_5 \\
    \begin{block}{c[rrr]}
    A & -1 & 0 & 1 \\ 
    B & 0 & 0 & 0 \\
    C & 1 & -1 & 0 \\
\end{block}
\end{blockarray} \hspace{.3cm} \text{and} \hspace{.2cm}
\begin{blockarray}{ccc}
    & R_2 & R_3 \\
    \begin{block}{c[rr]}
    A & -1 & 1 \\ 
    B & 1 & -1 \\
    C & 0 & 0 \\
\end{block}
\end{blockarray},
\end{equation*}
respectively. It can easily be verified that $s_1 = 2$ and $s_2 = 1$. Since we have equality $s=3=2+1=s_1+s_2$, the decomposition is independent.
\end{example}

Now, a method was developed in \cite{HDLC2021} to obtain an independent decomposition of a reaction network that employs a coordinate graph to find the desired decomposition. It was then proved that the obtained independent decomposition in \cite{HDLC2021} is already the finest (i.e.the independent decomposition with the maximum number of subnetworks) \cite{Hernandezetal2022}.

A key feature of independent decomposition is that it allows us to compute the positive steady states of the entire network by first deriving the positive steady states of each independent subnetwork and then intersecting them. This is formalized in the following theorem which is taken from \cite{FeinbergBook2019}.

\begin{theorem}
    Let $\mathcal{N}$ be a reaction network with kinetics $\mathcal{K}$ decomposed into subnetworks $\mathcal{N}_1,\mathcal{N}_2,\ldots,\mathcal{N}_\alpha$ and $\mathcal{K}_i$ be the restriction of $\mathcal{K}$ to reactions in $\mathcal{N}_i$. Then
    \begin{equation*}
        E_1 \cap E_2 \cap \cdots \cap E_\alpha \subseteq E
    \end{equation*}
    where $E$ is the set of positive steady states of the whole network while $E_i$ is the set of positive steady states of subnetwork $\mathcal{N}_i$. Furthermore, if the network decomposition is independent, then equality holds, i.e.,
    \begin{equation*}
        E_1 \cap E_2 \cap \cdots \cap E_\alpha = E.
    \end{equation*}
\end{theorem}

\subsection*{Network translation}

The parametrization method in \cite{JMP2019:parametrization} requires a CRN to be weakly reversible and deficiency zero. When these conditions are not satisfied, the reaction network must be transformed to obtain a network with the desired structural properties. This is done by the process of \textit{network translation}, first introduced in \cite{Johnston2014}, which modifies the graphical properties of the CRN but preserves the reaction vectors, and allows the original source complex to determine the kinetics so that the dynamics of the system is preserved.

Network translation can be visualized by the operation of adding or subtracting linear combinations of species, known as \textit{translation complexes}, from individual reactions \cite{JohnstonBurton2019}. Formally, we let $\alpha_k \in \mathbb{Z}_{\geq 0}^m$ to be a translation complex for reaction $k$. We represent the operation of translating the reaction $r_k=y_i\to y_j \in \mathcal{R}$ by the translation complex $\alpha_k$ as
\begin{equation}
    y_i \to y_j \hspace{1cm} (+\alpha_k) \label{eq:translate}
\end{equation}
for $k=1,\ldots,r$. Here, we add the translation complex to both reactant and product complexes to produce the translated reaction.

\begin{example}
\label{eq:translationscheme}
    Going back to Example \ref{main:example:CRN}, translating $R_2$ by $\alpha_2 = A+B$ and $R_3$ by $\alpha_3 = B+C$ gives us the following translated network: 
\begin{center}
\begin{tikzpicture}[baseline=(current  bounding  box.center)]
\tikzset{vertex/.style = {draw=none,fill=none}}
\tikzset{edge/.style = {->,> = latex', line width=0.10mm}} 
\node[vertex] (1) at  (-1.5,0) {$A+B$};
\node[vertex] (2) at  (1,0) {$B+C$};
\node[vertex] (3) at  (-0.2,-1.5) {$A+B+C$};
\node[vertex] (4) at  (2,0) {$A$};
\node[vertex] (5) at  (4,0) {$B$};

\draw [edge]  (1) to["$k_1$"] (2);
\draw [edge]  (3) to["$k_4$"] (1);
\draw [edge]  (2) to["$k_5$"] (3);
\draw [edge] ([yshift=2pt]4.east) -- ([yshift=2pt]5.west) node[midway, above] {$k_2$};
\draw [edge] ([yshift=-2pt]5.west) -- ([yshift=-2pt]4.east) node[midway, below] {$k_3$};
\end{tikzpicture}
\end{center}

Clearly, the translated network is weakly reversible. Furthermore, it can be verified that its deficiency is equal to zero, which is not the case for the original network.
\end{example}

The translation scheme in (\ref{eq:translate}) produces a new set of complexes, containing $y_i + \alpha_k$ and $y_j + \alpha_k$, which are called \textit{stoichiometric complexes} in the translated network. Since the original source complex determines the kinetics of the translated network, another set of complexes is produced called the \textit{kinetic complexes}. Therefore, the translated network can be thought of in terms of a generalized chemical reaction network (GCRN).

{
\subsection*{Generalized chemical reaction networks}
In this subsection, we formally define the concept of a GCRN, which is pioneered by M\"uller and Regensburger \cite{Muller2012}, and show how this is useful in deriving the analytic equilibria of a reaction network.

Let $G=(V,E)$ be a directed graph with vertex set $V$ and edge set $E \subset V \times V$. Furthermore, denote $V_s = \{ i | i \to j \in E\}$ to be the set of all source vertices in $G$. A GCRN is a directed graph $G=(V,E)$ with two maps
\begin{itemize}
    \item $y: V \to \mathbb{R}_{\geq 0}^m$ that assigns to each vertex a stoichiometric complex; and
    \item $\Tilde{y}: V_s \to \mathbb{R}_{\geq 0}^m$ that assigns to each vertex a kinetic complex.
\end{itemize}

Given a CRN with an associated graph $G$, we can construct a GCRN with graph $G'$ via network translation such that the resulting system of ODE matches that of the original CRN. In such case, the CRN and the GCRN are said to be \textit{dynamically equivalent}. The two maps $y$ and $\Tilde{y}$ give rise to two associated CRNs for a GCRN: the \textit{stoichiometric CRN} $(G',y)$ and the \textit{kinetic-order CRN} $(G',\Tilde{y})$. Consequently, we can define deficiencies for both the stoichiometric CRN and kinetic-order CRN, which are called \textit{effective deficiency} $\delta$ and \textit{kinetic deficiency} $\tilde{\delta}$, respectively. Both are calculated with the same formula $n-l-s$ but using the CRN associated with them.

There are two kinds of edges in a GCRN: phantom edge and effective edge. A \textit{phantom edge} connects identical stoichiometric complexes in the GCRN. Otherwise, it is an \textit{effective edge}. A phantom edge does not contribute to the system of ODEs since we would get $0$ as its corresponding reaction vector. Thus, we assign a dummy reaction rate constant $\sigma$, which is considered as a free parameter, for phantom edges. Lastly, we denote the sets of phantom edges and effective edges of the GCRN by $E^0$ and $E^*$, respectively. 

Given these, the main parametrization method is given in the theorem below \cite{HernandezPLOS2023,JMP2019:parametrization}.

\begin{theorem}
\label{thm:deficiency-based}
    Consider a weakly reversible translated network. Let $\mathcal{F}$ be any spanning forest containing all the nodes of the kinetic-order CRN. For each edge of $\mathcal{F}$, we define the kinetic difference as the vector produced by subtracting the head kinetic complex by the tail kinetic complex. Furthermore, $M$ is the matrix containing all the kinetic differences as rows where the entries per row are arranged according to the order of species. Let $H$ be a generalized inverse of $M$ (i.e., $MHM=M$). Finally, define $B$ and $C$ such that im$(B)=$ ker$(M)$ and ker$(B)=\{0\}$, and im$(C)=$ ker$\left(M^T\right)$ and ker$(C)=\{0\}$.
    \begin{itemize}
        \item (Case 1) If the kinetic deficiency is zero, the set of parametrized complex-balanced equilibria is given by
    \begin{equation*}
        \Bar{Z} = \left\{ \kappa\left( k^*,\sigma \right)^{H^T} \circ \tau^{B^T} | \sigma\in \mathbb{R}_{\geq 0}^{E^0}, \hspace{.2cm} \tau \in \mathbb{R}_{>0}^{m-\Tilde{s}}  \right\} \neq \emptyset
    \end{equation*}
        \item (Case 2) If the kinetic deficiency is positive, define the $\tilde{\delta}$ equations
        \begin{equation*}
            \kappa \left( k^*,k^0 \right)^C = 1^{\tilde{\delta} \times 1} 
        \end{equation*}
        and add them as additional conditions to the parametrization in Case 1.  
    \end{itemize}
     In this formalization, the product $\kappa\left( k^*,\sigma \right)^{H^T} \circ \tau^{B^T}$ is the Hadamard product (the number of components of $\kappa$ is the number of edges in a spanning forest, which can be effective of phantom) with the component of $\kappa$ associated with the edge $i\to i'$ as $\kappa_{i \to i'} = \dfrac{K_{i'}}{K_i}$ and tree constant $K_i$ as the sum (over all spanning trees of the kinetic-order CRN towards node $i$) of the products of the rate constants associated with the edges of each spanning tree. In addition, if the effective deficiency is zero, then the set of positive steady states of the original network is precisely $\Bar{Z}$.
\end{theorem}

\begin{example}
    In Example \ref{main:example:CRN} and its translation scheme in Example \ref{eq:translationscheme}, we define the GCRN as follows
\begin{center}
\begin{tikzpicture}[baseline=(current  bounding  box.center)]
\tikzset{vertex/.style = {draw=none,fill=none}}
\tikzset{edge/.style = {->,> = latex', line width=0.10mm}} 
\node[vertex] (1) at  (-1.5,0) {\shortstack{$A+B$\\$(A+B)$}};
\node[vertex] (2) at  (1.5,0) {\shortstack{$B$\\$(B+C)$}};
\node[vertex] (3) at  (0,-2) {\shortstack{$A$\\$(A+B+C)$}};

\node[vertex] (4) at  (3,0) {\shortstack{$A$\\$(A)$}};
\node[vertex] (5) at  (5,0) {\shortstack{$B$\\$(B)$}};

\draw [edge]  (1) to["$k_1$"] (2);
\draw [edge]  (2) to["$k_5$"] (3);
\draw [edge]  (3) to["$k_4$"] (1);
\draw [edge] ([yshift=2pt]4.east) -- ([yshift=2pt]5.west) node[midway, above] {$k_2$};
\draw [edge] ([yshift=-2pt]5.west) -- ([yshift=-2pt]4.east) node[midway, below] {$k_3$};
\end{tikzpicture}
\end{center}

Here, the complexes in parenthesis are the stoichiometric complexes while the complexes above them are the kinetic complexes. It can be verified that the GCRN has $\delta=0$ and $\tilde{\delta}=0$.
\end{example} 
}

\subsection*{Elementary flux modes}

Johnston and Burton \cite{JohnstonBurton2019} developed a computational method for performing structural translation through an elementary flux mode-based approach. An \textit{elementary flux mode} (EFM) represents a minimal, non-decomposable set of reactions that operate at steady state. In other words, EFMs are flux-balanced pathways that cannot be simplified in the sense that it is not possible to remove a subset of active reactions from an EFM and still be able to build a flux-balanced path using only the remaining active reactions \cite{KlapperSIAM2021}. 

Formally, a vector $e_i \in \mathbb{R}_{\geq 0}^r$ is an EFM of the CRN if $e_i \in \text{ker}(N)$ and $e_i$ is not a convex combination of any other $e_j,e_k \in \text{ker}(N) \cap \mathbb{R}_{\geq 0}^r$. The set $\text{ker}(N) \cap \mathbb{R}_{\geq 0}^r$ is the set of all admissible/feasible flux vectors, sometimes referred to as the \textit{elementary flux cone} $\text{cone}(\mathcal{E})$, where $\mathcal{E}$ is the set of all EFMs of a CRN. We say that $\mathcal{E}$ is \textit{unitary} if every entry in every $e_i \in \mathcal{E}$ is a one or a zero. Furthermore, $\mathcal{E}$ \textit{covers} the reaction set $\mathcal{R}$ is $\text{cone}(\mathcal{E}) \cap \mathbb{R}_{>0}^r \neq \emptyset$ \cite{JohnstonBurton2019}.

If the set of EFMs of a CRN is unitary and covers $\mathcal{R}$, then each EFM $e_i$ is completely determined by its support, $\text{supp}(e_i)$, which identifies the subset of reactions active in that mode. As a result, it is sufficient to represent an EFM by the positions of its nonzero entries, allowing us to use the binary vector notation (i.e., $e_i = (1,0,1,1,\ldots)$) and the corresponding reaction set notation (i.e., $e_i = \{ R_1,R_3,R_4,\ldots \}$) interchangeably.

The EFMs can be interpreted as the sets of reactions that, if taken in any order, would result in no net gain or loss of any species. In other words, each is a flux path through the reaction network with the property that the net total flux into and out of all nodes/vertices balances. Furthermore, EFMs are nondecomposable flux vectors from which all possible flux vectors of the network can be completely described \cite{KlapperSIAM2021}. That is, the EFMs are the extremal generators of $\text{cone}(\mathcal{E})$.

\begin{example}
\label{ex:efms}
    Note that the stoichiometric matrices of the network in Example \ref{main:example:CRN} and its weakly reversible, deficiency zero translation in Example \ref{eq:translationscheme} coincide since both networks have the same set of reaction vectors. Since EFMs are, by definition, elements of the null space of the stoichiometric matrix, the original and translated networks have the same set of EFMs. In particular, they have two EFMs, which are given by $e_1= \{ R_1,R_4,R_5 \}$ and $e_2 = \{ R_2,R_3 \}$.
\end{example}

\subsection*{Reaction-to-reaction graph}

A \textit{reaction-to-reaction graph} is a directed graph whose vertices correspond to the reactions of the CRN \cite{JohnstonBurton2019}. Directed cycles in the reaction-to-reaction graph are formed using the EFMs of the reaction network. Formally, a directed graph $G^\mathcal{R} = (V^\mathcal{R}, E^\mathcal{R})$ is a reaction-to-reaction graph of a CRN if $V^\mathcal{R}=\mathcal{R}$ and $E^\mathcal{R} \subset \mathcal{R} \times \mathcal{R}$. 

It would be convenient to use the mappings $s,p : \text{supp}(\mathcal{R}) \mapsto \text{supp}(\mathcal{C})$ in the subsequent definitions, where $s$ and $p$ map the source and product of each reaction to its corresponding complex, respectively. That is, we may use the notation $r_k = y_{s(k)}  \to y_{p(k)}$ instead of $r_k = y_i  \to y_j$. Now, we say that the CRN and the reaction-to-reaction graph are:
\begin{itemize}
    \item \textit{product-to-source compatible} (PS-compatible) if, for any $r_i=y_{s(i)} \to y_{p(i)}$ and $r_j=y_{s(j)} \to y_{p(j)}$, $(r_i,r_j) \in E^\mathcal{R}$ if and only if $y_{p(i)} = y_{s(j)}$. 
    \item \textit{common source compatible} (CS-compatible) if $y_{s(i)}=y_{s(j)}$ and $(r_k,r_i) \in E^\mathcal{R}$ implies $(r_k,r_j) \in E^\mathcal{R}$
    \item \textit{elementary flux mode compatible} (EM-compatible) if every EFM of the CRN corresponds to the vertices of a minimal directed cycle in the reaction-to-reaction graph, and vice versa
\end{itemize}

The theorem below relates the properties of PS-, CS-, and EM-compatibility to a weakly reversible, deficiency zero translation of a network. This is taken from \cite{JohnstonBurton2019}. 

\begin{theorem}
\label{thm: structuraltranslation}
    Consider a CRN $(\mathcal{S},\mathcal{C},\mathcal{R})$ with a set of elementary flux modes $\mathcal{E} = \{ e_1,e_2,\ldots,e_p \}$ which is unitary and covers $\mathcal{R}$. If there is a reaction-to-reaction graph $G^\mathcal{R} = (V^\mathcal{R},E^\mathcal{R})$ which is CS- and EM-compatible to $(\mathcal{S},\mathcal{C},\mathcal{R})$, then there is a CRN $(\mathcal{S},\mathcal{C}',\mathcal{R}')$ which is PS- CS, and EM-compatible with $G^\mathcal{R}$. Furthermore, $(\mathcal{S},\mathcal{C}',\mathcal{R}')$ is a weakly reversible, zero deficiency structural translation of $(\mathcal{S},\mathcal{C},\mathcal{R})$. In particular, the translation complexes $\alpha_i \in \mathbb{R}_{\geq 0}^m , i=1,\ldots,r,$ required to produce such a translation satisfy the following linear system, which is necessarily consistent:
    \begin{equation}
        \alpha_i - \alpha_j = y_{s(j)} - y_{p(i)}, \hspace{.5cm} (r_i,r_j) \in E^\mathcal{R}. \label{eq:2}
    \end{equation}
\end{theorem}

The goal then is to construct a reaction-to-reaction graph which are CS- and EM-compatible with a given CRN and then enforce PS-compatibility to construct a weakly reversible, deficiency zero translation via solving a linear system. 

\begin{example}
Consider the reaction-to-reaction graph below.
\begin{center}
\begin{tikzpicture}[baseline=(current  bounding  box.center)]
\tikzset{vertex/.style = {draw=none,fill=none}}
\tikzset{edge/.style = {->,> = latex', line width=0.10mm}} 
\node[vertex] (1) at  (-1.5,0) {$R_1$};
\node[vertex] (2) at  (1.5,0) {$R_5$};
\node[vertex] (3) at  (0,-1) {$R_4$};

\node[vertex] (4) at  (3,0) {$R_2$};
\node[vertex] (5) at  (5,0) {$R_3$};

\draw [edge]  (1) to (2);
\draw [edge]  (2) to (3);
\draw [edge]  (3) to (1);
\draw [edge] ([yshift=2pt]4.east) -- ([yshift=2pt]5.west) node[midway, above] {};
\draw [edge] ([yshift=-2pt]5.west) -- ([yshift=-2pt]4.east) node[midway, below] {};
\end{tikzpicture}
\end{center}

Indeed, the EFMs computed in Example \ref{ex:efms} given by $e_1 = \{ R_1,R_4,R_5 \}$ and $e_2 = \{ R_2,R_3 \}$ correspond to minimal cycles in the reaction-to-reaction graph. Moreover, since there are no reactions with common source complex in the network in Example \ref{main:example:CRN}, then the CRN and the reaction-to-reaction graph are automatically CS-compatible. However, the reaction-to-reaction graph is not PS-compatible to the CRN. For instance, consider $(R_1,R_2) \in E^\mathcal{R}$. We have $y_{p(1)} = C+B \neq A=y_{s(2)}$. Nevertheless, the reaction-to-reaction graph is EM- and CS-compatible with the network in Example \ref{main:example:CRN}. Moreover, the set of EFMs $\mathcal{E} = \{ e_1,e_2 \}$ is unitary and covers $\mathcal{R}$. By Theorem, \ref{thm: structuraltranslation}, there exists a weakly reversible, deficiency zero translation of the network. Indeed, such a translation was constructed in Example \ref{eq:translationscheme}.
\end{example}


\section*{Acknowledgments}
The authors acknowledge the Office of the Chancellor of the University of the Philippines Diliman, through the Office of the Vice Chancellor for Research and Development, for funding support through the Outright Research Grant (242414 ORG).

\section*{Author contributions}
{\bf Conceptualization:} Exequiel Jun V. Villejo, Aurelio A. de los Reyes V, Bryan S. Hernandez.\\ 
\noindent{\bf Funding acquisition:} Bryan S. Hernandez.\\
\noindent{\bf Methodology:} Exequiel Jun V. Villejo, Bryan S. Hernandez.\\
\noindent{\bf Project administration:} Bryan S. Hernandez.\\
\noindent{\bf Software:} Exequiel Jun V. Villejo.\\
\noindent{\bf Supervision:} Bryan S. Hernandez.\\
\noindent{\bf Visualization:} Exequiel Jun V. Villejo, Aurelio A. de los Reyes V, Bryan S. Hernandez.\\
\noindent{\bf Writing – original draft:} Exequiel Jun V. Villejo, Aurelio A. de los Reyes V, Bryan S. Hernandez.\\
\noindent{\bf Writing – review \& editing:} Exequiel Jun V. Villejo, Aurelio A. de los Reyes V, Bryan S. Hernandez.


%
%
%


\section*{Supporting information}


\paragraph*{S1 Text.}
Supplementary Methods, Table, and Figures.

\end{document}